\font\subtit=cmr12

\font\name=cmr8
%%%%%%%%%%%%%%%%%%  tex macros for preprints, cm version %%%%%%%%%%%%%%
%                     (P. Ginsparg, last updated 9/91)
%                if confused, type `b' in response to query 
%
%---------------------------------------------------------------------%
%% site dependent options: 
%% \unredoffs and \redoffs define horizontal and vertical offsets 
%% respectively for unreduced and reduced modes. \speclscape defines
%% the \special{} call that sets printer to landscape (sideways) mode.
%% from standard set below, leave uncommented as appropriate or redefine
%
%%% next 400dpi
%\def\unredoffs{} \def\redoffs{\voffset=-.31truein\hoffset=-.48truein}
%\def\speclscape{\special{landscape}}
%
%%% apple lw
\def\unredoffs{} \def\redoffs{\voffset=-.31truein\hoffset=-.59truein}
\def\speclscape{\special{ps: landscape}}
%
%%% qms lasergrafix:
%\def\unredoffs{} \def\redoffs{\voffset=-.4truein\hoffset=.125truein}
%\def\speclscape{\special{qms: landscape}}
%
%%% saclay A4 paper:
%\def\unredoffs{\hoffset-.14truein\voffset-.2truein} 
%\def\redoffs{\voffset=-.45truein\hoffset=-.21truein} 
%\def\speclscape{\special{landscape}}
%
%---------------------------------------------------------------------%
%
\newbox\leftpage \newdimen\fullhsize \newdimen\hstitle \newdimen\hsbody
\tolerance=1000\hfuzz=2pt
\catcode`\@=11 % This allows us to modify PLAIN macros.
\def\bigans{b }
%\message{ big or little (b/l)? }\read-1 to\answ
%
\def\answ{b }
\ifx\answ\bigans\message{(This will come out unreduced.}
%originariamnete \baselineskip=16pt
\magnification=1200\unredoffs\baselineskip=24pt plus 2pt minus 1pt
\hsbody=\hsize \hstitle=\hsize %take default values for unreduced format
\else\message{(This will be reduced.} \let\l@r=L
\magnification=1000\baselineskip=16pt plus 2pt minus 1pt \vsize=7truein
\redoffs \hstitle=8truein\hsbody=4.75truein\fullhsize=10truein\hsize=\hsbody
\output={\ifnum\pageno=0 %%% This is the HUTP version
  \shipout\vbox{\speclscape{\hsize\fullhsize\makeheadline}
    \hbox to \fullhsize{\hfill\pagebody\hfill}}\advancepageno
  \else
  \almostshipout{\leftline{\vbox{\pagebody\makefootline}}}\advancepageno 
  \fi}
\def\almostshipout#1{\if L\l@r \count1=1 \message{[\the\count0.\the\count1]}
      \global\setbox\leftpage=#1 \global\let\l@r=R
 \else \count1=2
  \shipout\vbox{\speclscape{\hsize\fullhsize\makeheadline}
      \hbox to\fullhsize{\box\leftpage\hfil#1}}  \global\let\l@r=L\fi}
\fi
%---------------------------------------------------------------------
%
\newcount\yearltd\yearltd=\year\advance\yearltd by -1900

\def\Date#1{\vfill\leftline{#1}\tenpoint\supereject\global\hsize=\hsbody%
\footline={\hss\tenrm\folio\hss}}% 	restores pagenumbers
%
%       use following instead of \Date on the preliminary draft, 
%       puts date/time on each page in big mode, writes labels in margins

\def\draftmode{\message{ DRAFTMODE }\def\draftdate{{\rm preliminary draft:
\number\month/\number\day/\number\yearltd\ \ \hourmin}}%
\headline={\hfil\draftdate}\writelabels\baselineskip=20pt plus 2pt minus 2pt
 {\count255=\time\divide\count255 by 60 \xdef\hourmin{\number\count255}
  \multiply\count255 by-60\advance\count255 by\time
  \xdef\hourmin{\hourmin:\ifnum\count255<10 0\fi\the\count255}}}
%       use \nolabels to get rid of eqn, ref, and fig labels in draft mode
\def\nolabels{\def\wrlabeL##1{}\def\eqlabeL##1{}\def\reflabeL##1{}}
\def\writelabels{\def\wrlabeL##1{\leavevmode\vadjust{\rlap{\smash%
{\line{{\escapechar=` \hfill\rlap{\sevenrm\hskip.03in\string##1}}}}}}}%
\def\eqlabeL##1{{\escapechar-1\rlap{\sevenrm\hskip.05in\string##1}}}%
\def\reflabeL##1{\noexpand\llap{\noexpand\sevenrm\string\string\string##1}}}
\nolabels
%
% tagged sec numbers
\global\newcount\secno \global\secno=0
\global\newcount\meqno \global\meqno=1
\def\newsec#1{\global\advance\secno by1\message{(\the\secno. #1)}
%\ifx\answ\bigans \vfill\eject \else \bigbreak\bigskip \fi  %if desired
\global\subsecno=0\eqnres@t\noindent{\bf\uppercase\expandafter{\romannumeral\the\secno}. #1}
\writetoca{{\secsym}{#1}}\par\nobreak\medskip\nobreak}
\def\eqnres@t{\xdef\secsym{\the\secno.}\global\meqno=1\bigbreak\bigskip}
\def\sequentialequations{\def\eqnres@t{\bigbreak}}\xdef\secsym{}
\global\newcount\subsecno \global\subsecno=0
\def\subsec#1{\global\advance\subsecno by1\message{(\secsym\the\subsecno. #1)}
\ifnum\lastpenalty>9000\else\bigbreak\fi
\noindent{\it\secsym\the\subsecno. #1}\writetoca{\string\quad 
{\secsym\the\subsecno.} {#1}}\par\nobreak\medskip\nobreak}
\def\appendix#1#2{\global\meqno=1\global\subsecno=0\xdef\secsym{\hbox{#1.}}
\bigbreak\bigskip\noindent{\bf Appendix #1. #2}\message{(#1. #2)}
\writetoca{Appendix {#1.} {#2}}\par\nobreak\medskip\nobreak}
%
%       \eqn\label{a+b=c}	gives displayed equation, numbered
%				consecutively within sections.
%     \eqnn and \eqna define labels in advance (of eqalign?)
%
\def\eqnn#1{\xdef #1{(\secsym\the\meqno)}\writedef{#1\leftbracket#1}%
\global\advance\meqno by1\wrlabeL#1}
\def\eqna#1{\xdef #1##1{\hbox{$(\secsym\the\meqno##1)$}}
\writedef{#1\numbersign1\leftbracket#1{\numbersign1}}%
\global\advance\meqno by1\wrlabeL{#1$\{\}$}}
\def\eqn#1#2{\xdef #1{(\secsym\the\meqno)}\writedef{#1\leftbracket#1}%
\global\advance\meqno by1$$#2\eqno#1\eqlabeL#1$$}
%
%			 footnotes
\newskip\footskip\footskip14pt plus 1pt minus 1pt %sets footnote baselineskip
\def\footnotefont{\ninepoint}\def\f@t#1{\footnotefont #1\@foot}
\def\f@@t{\baselineskip\footskip\bgroup\footnotefont\aftergroup\@foot\let\next}
\setbox\strutbox=\hbox{\vrule height9.5pt depth4.5pt width0pt}
\global\newcount\ftno \global\ftno=0
\def\foot{\global\advance\ftno by1\footnote{$^{\the\ftno}$}}
%
%say \footend to put footnotes at end
%will cause problems if \ref used inside \foot, instead use \nref before
\newwrite\ftfile   
\def\footend{\def\foot{\global\advance\ftno by1\chardef\wfile=\ftfile
$^{\the\ftno}$\ifnum\ftno=1\immediate\openout\ftfile=foots.tmp\fi%
\immediate\write\ftfile{\noexpand\smallskip%
\noexpand\item{f\the\ftno:\ }\pctsign}\findarg}%
\def\footatend{\vfill\eject\immediate\closeout\ftfile{\parindent=20pt
\centerline{\bf Footnotes}\nobreak\bigskip\input foots.tmp }}}
\def\footatend{}

%
%     \ref\label{text}
% generates a number, assigns it to \label, generates an entry.
% To list the refs on a separate page,  \listrefs
%
\global\newcount\refno \global\refno=1
\newwrite\rfile
\def\ref{[\the\refno]\nref}
\def\nref#1{\xdef#1{[\the\refno]}\writedef{#1\leftbracket#1}%
\ifnum\refno=1\immediate\openout\rfile=refs.tmp\fi
\global\advance\refno by1\chardef\wfile=\rfile\immediate
\write\rfile{\noexpand\item{#1\ }\reflabeL{#1\hskip.31in}\pctsign}\findarg}
%	horrible hack to sidestep tex \write limitation
\def\findarg#1#{\begingroup\obeylines\newlinechar=`\^^M\pass@rg}
{\obeylines\gdef\pass@rg#1{\writ@line\relax #1^^M\hbox{}^^M}%
\gdef\writ@line#1^^M{\expandafter\toks0\expandafter{\striprel@x #1}%
\edef\next{\the\toks0}\ifx\next\em@rk\let\next=\endgroup\else\ifx\next\empty%
\else\immediate\write\wfile{\the\toks0}\fi\let\next=\writ@line\fi\next\relax}}
\def\striprel@x#1{} \def\em@rk{\hbox{}} 
\def\lref{\begingroup\obeylines\lr@f}
\def\lr@f#1#2{\gdef#1{\ref#1{#2}}\endgroup\unskip}
\def\semi{;\hfil\break}
\def\addref#1{\immediate\write\rfile{\noexpand\item{}#1}} %now unnecessary
\def\startrefs#1{\immediate\openout\rfile=refs.tmp\refno=#1}
\def\xref{\expandafter\xr@f}\def\xr@f[#1]{#1}
\def\refs#1{\count255=1[\r@fs #1{\hbox{}}]}
\def\r@fs#1{\ifx\und@fined#1\message{reflabel \string#1 is undefined.}%
\nref#1{need to supply reference \string#1.}\fi%
\vphantom{\hphantom{#1}}\edef\next{#1}\ifx\next\em@rk\def\next{}%
\else\ifx\next#1\ifodd\count255\relax\xref#1\count255=0\fi%
\else#1\count255=1\fi\let\next=\r@fs\fi\next}
%

%
% this is ugly, but moore insists
\newwrite\ffile\global\newcount\figno \global\figno=1
\def\fig{fig.~\the\figno\nfig}
\def\nfig#1{\xdef#1{fig.~\the\figno}%
\writedef{#1\leftbracket fig.\noexpand~\the\figno}%
\ifnum\figno=1\immediate\openout\ffile=figs.tmp\fi\chardef\wfile=\ffile%
\immediate\write\ffile{\noexpand\medskip\noexpand\item{Fig.\ \the\figno. }
\reflabeL{#1\hskip.55in}\pctsign}\global\advance\figno by1\findarg}
\def\xfig{\expandafter\xf@g}\def\xf@g fig.\penalty\@M\ {}
\def\figs#1{figs.~\f@gs #1{\hbox{}}}
\def\f@gs#1{\edef\next{#1}\ifx\next\em@rk\def\next{}\else
\ifx\next#1\xfig #1\else#1\fi\let\next=\f@gs\fi\next}
\newwrite\lfile
{\escapechar-1\xdef\pctsign{\string\%}\xdef\leftbracket{\string\{}
\xdef\rightbracket{\string\}}\xdef\numbersign{\string\#}}

\def\writestop{\def\writestoppt{\immediate\write\lfile{\string\pageno%
\the\pageno\string\startrefs\leftbracket\the\refno\rightbracket%
\string\def\string\secsym\leftbracket\secsym\rightbracket%
\string\secno\the\secno\string\meqno\the\meqno}\immediate\closeout\lfile}}
\def\writestoppt{}\def\writedef#1{}
\def\seclab#1{\xdef #1{\the\secno}\writedef{#1\leftbracket#1}\wrlabeL{#1=#1}}
\def\subseclab#1{\xdef #1{\secsym\the\subsecno}%
\writedef{#1\leftbracket#1}\wrlabeL{#1=#1}}
\newwrite\tfile \def\writetoca#1{}
\def\leaderfill{\leaders\hbox to 1em{\hss.\hss}\hfill}
%	use this to write file with table of contents
\def\writetoc{\immediate\openout\tfile=toc.tmp 
   \def\writetoca##1{{\edef\next{\write\tfile{\noindent ##1 
   \string\leaderfill {\noexpand\number\pageno} \par}}\next}}}
%       and this lists table of contents on second pass

%
\catcode`\@=12 % at signs are no longer letters
%
%	Unpleasantness in calling in abstract and title fonts
\edef\tfontsize{\ifx\answ\bigans scaled\magstep3\else scaled\magstep4\fi}
 \tfontsize  \tfontsize
 \tfontsize \font\titlei=cmmi10 \tfontsize
\font\titleis=cmmi7 \tfontsize \font\titleiss=cmmi5 \tfontsize
\font\titlesy=cmsy10 \tfontsize \font\titlesys=cmsy7 \tfontsize
\font\titlesyss=cmsy5 \tfontsize  \tfontsize
\skewchar\titlei='177 \skewchar\titleis='177 \skewchar\titleiss='177
\skewchar\titlesy='60 \skewchar\titlesys='60 \skewchar\titlesyss='60
 \ifx\answ\bigans\else scaled\magstep1\fi
\ifx\answ\bigans\def\abstractfont{\tenpoint}\else
\font\abssl=cmsl10 scaled \magstep1
\font\absrm=cmr10 scaled\magstep1 \font\absrms=cmr7 scaled\magstep1
\font\absrmss=cmr5 scaled\magstep1 \font\absi=cmmi10 scaled\magstep1
\font\absis=cmmi7 scaled\magstep1 \font\absiss=cmmi5 scaled\magstep1
\font\abssy=cmsy10 scaled\magstep1 \font\abssys=cmsy7 scaled\magstep1
\font\abssyss=cmsy5 scaled\magstep1 \font\absbf=cmbx10 scaled\magstep1
\skewchar\absi='177 \skewchar\absis='177 \skewchar\absiss='177
\skewchar\abssy='60 \skewchar\abssys='60 \skewchar\abssyss='60
\def\abstractfont{\def\rm{\fam0\absrm}% switch to abstract font
\textfont0=\absrm \scriptfont0=\absrms \scriptscriptfont0=\absrmss
\textfont1=\absi \scriptfont1=\absis \scriptscriptfont1=\absiss
\textfont2=\abssy \scriptfont2=\abssys \scriptscriptfont2=\abssyss
\textfont\itfam=\bigit \def\it{\fam\itfam\bigit}\def\footnotefont{\tenpoint}%
\textfont\slfam=\abssl \def\sl{\fam\slfam\abssl}%
\textfont\bffam=\absbf \def\bf{\fam\bffam\absbf}\rm}\fi
\def\tenpoint{\def\rm{\fam0\tenrm}% switch back to 10-point type
\textfont0=\tenrm \scriptfont0=\sevenrm \scriptscriptfont0=\fiverm
\textfont1=\teni  \scriptfont1=\seveni  \scriptscriptfont1=\fivei
\textfont2=\tensy \scriptfont2=\sevensy \scriptscriptfont2=\fivesy
\textfont\itfam=\tenit \def\it{\fam\itfam\tenit}\def\footnotefont{\ninepoint}%
\textfont\bffam=\tenbf \def\bf{\fam\bffam\tenbf}\def\sl{\fam\slfam\tensl}\rm}
\font\ninerm=cmr9 \font\sixrm=cmr6 \font\ninei=cmmi9 \font\sixi=cmmi6 
\font\ninesy=cmsy9 \font\sixsy=cmsy6 \font\ninebf=cmbx9 
\font\nineit=cmti9 \font\ninesl=cmsl9 \skewchar\ninei='177
\skewchar\sixi='177 \skewchar\ninesy='60 \skewchar\sixsy='60 
\def\ninepoint{\def\rm{\fam0\ninerm}% switch to footnote font
\textfont0=\ninerm \scriptfont0=\sixrm \scriptscriptfont0=\fiverm
\textfont1=\ninei \scriptfont1=\sixi \scriptscriptfont1=\fivei
\textfont2=\ninesy \scriptfont2=\sixsy \scriptscriptfont2=\fivesy
\textfont\itfam=\ninei \def\it{\fam\itfam\nineit}\def\sl{\fam\slfam\ninesl}%
\textfont\bffam=\ninebf \def\bf{\fam\bffam\ninebf}\rm} 
%
%---------------------------------------------------------------------
%

\hyphenation{anom-aly anom-alies coun-ter-term coun-ter-terms}
\def\inv{^{\raise.15ex\hbox{${\scriptscriptstyle -}$}\kern-.05em 1}}

\def\Dsl{\,\raise.15ex\hbox{/}\mkern-13.5mu D} %this one can be subscripted
\def\dsl{\raise.15ex\hbox{/}\kern-.57em\partial}

\font\bigit=cmti10 scaled \magstep1
 %pound sterling
\def\lspace{\ifx\answ\bigans{}\else\qquad\fi}
\def\lbspace{\ifx\answ\bigans{}\else\hskip-.2in\fi} % $$\lbspace...$$
\def\boxeqn#1{\vcenter{\vbox{\hrule\hbox{\vrule\kern3pt\vbox{\kern3pt
	\hbox{${\displaystyle #1}$}\kern3pt}\kern3pt\vrule}\hrule}}}
\def\mbox#1#2{\vcenter{\hrule \hbox{\vrule height#2in
		\kern#1in \vrule} \hrule}}  %e.g. \mbox{.1}{.1}
%	matters of taste
%\def\tilde{\widetilde} \def\bar{\overline} \def\hat{\widehat}
%
% some sample definitions
  %     curly letters

\def\darr#1{\raise1.5ex\hbox{$\leftrightarrow$}\mkern-16.5mu #1}
 %pound sterling

 %puts a small half in a displayed eqn
\def\roughly#1{\raise.3ex\hbox{$#1$\kern-.75em\lower1ex\hbox{$\sim$}}}

\def\plb#1#2#3#4{#1, {\it Phys. Lett.} {\bf {#2}}B (#3), #4}
\def\npb#1#2#3#4{#1, {\it Nucl. Phys.} {\bf B{#2}} (#3), #4}
\def\npbib#1#2#3{{\it Nucl. Phys.} {\bf B{#1}} (#2), #3}
\def\prl#1#2#3#4{#1, {\it Phys. Rev. Lett.} {\bf {#2}} (#3), #4}

\def\cmp#1#2#3#4{#1, {\it Comm. Math. Phys.} {\bf {#2}} (#3), #4}
\def\lmp#1#2#3#4{#1, {\it Lett. Math. Phys.} {\bf {#2}} (#3), #4}

\def\mpl#1#2#3#4{#1, {\it Mod. Phys. Lett.} {\bf A{#2}} (#3), #4}
\def\ijmpa#1#2#3#4{#1, {\it Int. Jour. Mod. Phys.} {\bf A{#2}} (#3), #4}

\def\jgp#1#2#3#4{#1, {\it Journal Geom. Phys.} {\bf {#2}} (#3), #4}

%corrections get here
\null\vskip 1truecm
\centerline{\subtit
GENERALIZED WEIERSTRASS KERNELS ON THE}
\centerline{\subtit INTERSECTION OF TWO COMPLEX
HYPERSURFACES}

\vskip 1truecm
\centerline{F{\name RANCO} F{\name ERRARI}}
\smallskip \centerline{\it Center for Theoretical Physics, Laboratory for
Nuclear Science and}
\smallskip \centerline{\it Department of Physics, Massachusetts Institute of
Technology}
\smallskip \centerline{\it Cambridge, Massachusetts 02139, USA}
\smallskip \centerline{\it Physics Department,
University of Szczecin, Wielkopolska 15, 70-451 Szczecin, Poland\foot{E-mail:
fferrari@univ.szczecin.pl}.}
\smallskip \centerline{\it Physics Department,
University of Trento and INFN, Gruppo Coll. di Trento, 38050 Povo, Italy.}
\smallskip
\vskip 1cm
\centerline{ABSTRACT}
{\narrower \abstractfont
On plane algebraic curves the so-called
Weierstrass kernel plays the same role of the
Cauchy kernel on the complex plane. A straightforward prescription
to construct the Weierstrass kernel is known since one century.
How can it be extended to the case of more general curves obtained from
the intersection of hypersurfaces in a $n$ dimensional complex space?
This problem is solved in this work in the case $n=3$.
As an application, the correlation functions of bosonic string theories
are constructed on a canonical curve of genus four.
} 
%\Date{December 1996}
\Date{December 2000}
%\pageno=0
%\draft
\vfill\eject
\newsec {INTRODUCTION}
\vskip 1cm
Algebraic curves  and $n-$sheeted branched covers on the projective
line {\bf P}$^1$
provide an explicit
representation of abstract Riemann surfaces.
Besides being an active
field of research
in mathematics 
\ref\fakra{H. Farkas and I. Kra, {\it Riemann Surfaces},
Springer Verlag, 1980.}\nref\fay{J.  D.  Fay, Theta Functions on Riemann
Surfaces, Lecture Notes in Mathematical Physics no.  352, Springer Verlag,
1973.}\nref\zvero{E.I. Zverovich, {\it Russ.
Math. Surv.} {\bf 21} 1971, 99.}\nref\cenr{F.  Enriques and O.
Chisini, {\it Lezioni sulla Teoria
Geometrica delle Equazioni e delle Funzioni Algebriche}, Zanichelli,
Bologna (in italian).}\nref\hareis{J. Harris and
D. Eisenbud, {\it Bull. Am. Math. Soc.}
{\bf 21} (2) (1989), 205; B. Dubrovin,
{\it Jour. Diff. Geom.} Suppl. 4 (1998), 213; E. Ballico,
C. Keem, {\it Israel Jour.
Math.} {\bf 104} (1998), 355; C. Ciliberto, J. Harris,
{\it Commun. Algebra} {\bf 27} (3) (1999), 2197.}--\ref\grha{
P. Griffiths
and J. Harris, {\it Principles of Algebraic
Geometry}, John Wiley \& Sons, New York 1978.}, since more than a decade
they have been also successfully applied
in several different topics of theoretical
and mathematical physics
\ref\hyel{Some Applications
of hyperelliptic curves in perturbative string theory:
R. Iengo and C. J. Zhu, {\it Jour. High Energy Phys.}  {\bf 4} (2000),
U651; R. Iengo and C.-J. Zhu, JHEP {\bf 9906} (1999), 011.
\npb{D.  Lebedev and A.  Morozov}{302}{1986}{163}; \plb
{E. Gava, R. Iengo and G. Sotkov}{207}{1988}{283}; \plb{A. Yu. Morozov and
A. Perelomov} {197}{1987} {115}; \npb{D.  Montano}{297}{1988}{125};
F. Ferrari, {\it Fizika} {\bf 21} (1989), 32; \mpl{J. Sobczyk}{6}{1991}{1103};
\npb{D.  Montano}{297}{1988}{125}; \npb{E. Gava,
R.  Iengo and C.  J.  Zhu}{323}{1989}{585};
\plb{R.  Jengo and C.  J.  Zhu}{212}{1988}{313};
V. G. Knizhnik,
{\it Phys. Lett.} {\bf 196B} (1987), 473;
D. J. Gross and P. F. Mende, {\it Nucl. Phys.} {\bf B303} (1988), 407.
}\nref\nptasp{Some applications of algebraic curves to various
non-perturbative aspects
of string theory and related problems:
W. Lerche, {\it  On the Heterotic Theory Duality in Eight-Dimensions},
Proceedings of the TMR Summer School on Progress in String Theory and
M-Theory (Cargese 99), Cargese, Corsica, France, 24 May - 5 Jun 1999,
hep-th/9910207;
W. Lerche and S. Stieberger, {\it Adv. Theor. Math. Phys.} {\bf 2} (1998),
1105, hep-th/9804176; W. Lerche,
S. Stieberger
and N. P. Warner, {\it Quartic Gauge Couplings from K3 Geometry},
Preprint CERN-TH/98-378, hep-th/9811228; 
%\cite{Bertoldi:2000hj}:
%\bibitem{Bertoldi:2000hj}
G.~Bertoldi, J.~M.~Isidro, M.~Matone and P.~Pasti,
%``The concept of a noncommutative Riemann surface,''
Phys.\ Lett.\ {\bf B484} (2000) 323
[hep-th/0003200];
%%CITATION = HEP-TH 0003200;%%
%\cite{Bertoldi:2000ua}:
%\bibitem{Bertoldi:2000ua}
G.~Bertoldi, J.~M.~Isidro, M.~Matone and P.~Pasti,
%``Noncommutative Riemann surfaces,''
hep-th/0003131;
%%CITATION = HEP-TH 0003131;%%
%\cite{Matone:1995tj}:
%\bibitem{Matone:1995tj}
M.~Matone,
%``Uniformization theory and 2-D gravity. 1. Liouville action and
%intersection numbers,''
Int.\ J.\ Mod.\ Phys.\ {\bf A10} (1995) 289
[hep-th/9306150];
%%CITATION = HEP-TH 9306150;%%
\npb{N. Seiberg and E. Witten}{430}{1994}{485};
\npb{N. Seiberg and E. Witten}{431}{1994}{484};
\jgp{J. M. Isidro}{29}{1999}{334};
\cmp{E. Witten}{113}{1988}{529}.}\nref\cfft{
Some applications of algebraic curves in conformal field theories:
\npb{Al. B. Zamolodchikov}{285}{1987}{481};
\ijmpa{M. A.  Bershadsky and A.  O.  Radul}{2}{1987}{165};
\cmp{V. G. Knizhnik} {112}{1987} {587};
L. Borisov, M.B. Halpern and C. Schweigert, {\it Int. Jour. Mod. Phys.}
{\bf A13} (1) (1998), 125;
S. A. Apikian and C. J. Efthimiou, {\it Int. Jour.
Mod. Phys.}, {\bf A12} 1997, 4291, hep-th/9610051;
\plb{S. A. Apikian and C. J. Efthimiou}{383}{1996}(397);
C. J. Efthimiou and D. A. Spector, {\it A Collection of Exercises in
Two-Dimensional Physics, Part 1},
Preprint CLNS-99-1612, hep-th/0003190; \plb{C. Crnkovic, G. M. Sotkov and M.
Stanishkov}{220}{1989}{397};
M. A. Bershadsky and A. Radul,
{\it Phys. Lett.} {\bf 193B} (1987), 21.}
\nref\kholvil{Some applications of algebraic curves to polymer
physics: A. L. Kholodenko
and T. A. Vilgis, {\it Phys. Rep.} {\bf 298} (1998), 251;
S. Nechaev,
{\it Int. Jour. Mod. Phys.} {\bf B4} (1990), 1809;
{\it Statistics of Knots and Entangled Random Walks},
extended version of lectures presented at Les Houches
1998 Summer School on {\it Topological Aspects of Low
Dimensional Systems}, July 7 - 31, 1998,
cond-mat/9812205.}--\ref\others{Some other applications of algebraic curves:
\lmp{L. A. Takhtajan}{52}{2000}{79};
\npb{I. Bakas and K. Sfetsos}{573}{2000}{768}; S. M. Sergeev, {\it Jour.
Nonlinear. Math. Phy.} {\bf 7} (2000), 57; B. Crespi, S. J. Chang, K. J. Shi,
{\it Jour. Math. Phys.} {\bf 34} (6) (1993), 2257;
\cmp{S. Cordes, G. Moore and S. Rangoolam}{185}{1997}{543};
\npb{L. Dixon, D. Friedan, E. Martinec and  S Shenker}{282}{1987}{13};
V. Marotta and A. Sciarrino, {\it Mod. Phys. Lett.} {\bf A13} (1998), 2863;
V. Marotta, {\it Nucl.Phys.} {\bf B527} (1998), 717-737, hep-th/9702143.}.

All the above physical applications involve only
plane affine algebraic curves. Moreover,
most of the works are further restricted
to the particular case of
hyperelliptic curves.
The reason is that the latter are very well known in mathematical
literature. 
An invaluable source of hyperelliptic formulas can be found for
instance in the original Weierstrass lecture notes
\ref\weier{K. Weierstrass, Vorlesungen \"uber di Theorie der Abelschen
Transcendenten, Math. Werke. Vol. 4, Berlin 1902.}.
As an example, we give the prescription of Weierstrass for
differentials of the second kind.
Already in the case of the so-called $Z_n$ symmetric curves,
which are the most straightforward generalization of hyperelliptic curves,
superstring calculations like those of \hyel\ are impossible, since
it has not yet been understood
how to construct sections with half-integer
spins or chiral determinants
\ref\fms{\npb{D. Friedan, E. Martinec and S.
Shenker}{271}{1986}{93}.}--\ref\veal{\npb{E. Verlinde
and H. Verlinde}{288}{1987}
{357}; \npb{L. Alvarez-Gaum\'e, C. Gomez, P. Nelson, G. Sierra and C. Vafa}
{311}{1988}{333}.}
of free fermions with given
spin structures. An exception is provided by more exotic
fermions with $1/n$ characteristics, which have been solved in \ref\brss{
\cmp{M. Bershadsky and M. Radul}{116}{1988}{689}.}.
Also neglecting fermions and sticking to bosonic string theory,
the computation of the partition function remains a very complicated task
on non-hyperelliptic algebraic curves.

The explicitness of algebraic curves is indeed of great
advantage in understanding the physical aspects of
theories like superstrings whenever Riemann surfaces are involved,
but, on the other side, it is also the main responsible
of the above mentioned difficulties. Using different approaches, like for
instance theta functions, the problems are somewhat hidden, but they reappear
in other forms, in particular when one needs to explore a limited portion
of the moduli space of Riemann surfaces such as that spanned by
$Z_n$ symmetric curves.

A considerable
effort to understand conformal field theories on
plane
algebraic curves and to find new applications of non-hyperelliptic curves
has been made by the author over the
past ten years, partly in collaboration with J. Sobczyk and W. Urbanik.
In particular, the applications of branched covers of
{\bf P}$^1$ have been investigated,
in which the curves are projected on the complex projective line.
Every compact Riemann surface can be represented in this way.
Branched coverings are also closely related to affine curves
in the two dimensional complex plane {\bf C}$^2$.
The bosonic string theories on a general non-hyperelliptic
algebraic curve of genus three
have been treated in refs. \ref\ffstearI{
F. Ferrari,
{\it Int. Jour. Mod. Phys.} {\bf A5} (1990), 2799.}--\ref\ffstearII{
F. Ferrari,
{\it Int. Jour. Mod. Phys.} {\bf A7} (1992), 5131.}, computing the correlation
functions of the theory and its
partition function, the latter up to a theta constant.
Moreover, the partition function 
has been exactly derived
on $Z_3$ symmetric curves.
in \ffstearII.
Even in this case, where
the chiral
determinants are well known
\ref\bbrr{\ijmpa{M. A.  Bershadsky and A.  O.  Radul}{2}{1987}{165}.},
the calculation of the partition function is not simple.
Indeed, one has still to use the Rauch's variational method
\ref\rauch{H. Rauch, {\it Comm. Pure Appl. Math.} {\bf 12} (1959), 543.}
combined with the Beilinson-Manin formula \ref\beiman{\cmp{A. A. Beilinson
and Yu. I. Manin}{107}{1986}{359}.}
in order to determine the
explicit form of the period matrix.
The same variational approach of Rauch has been exploited
in
\ref\naka{A. Nakayashiki, {\it Publ. Res. Inst. Math. Sci.} {\bf 34}
(1998), 439; {\it ibid.} {\bf 33} (1997), 987.} to rederive in an elegant
way the Thomae formulae of \bbrr\ on $Z_n$ symmetric curves.

Despite the difficulties of evaluating the partition functions
of conformal field theories,
algebraic curves become very convenient in the construction of
meromorphic tensors with poles and zeros of given order at given points,
like for example the correlation functions of the free
fields appearing in the bosonic string action \fms.
On general plane algebraic curves, this problem has been solved
for the $b-c$ systems in
\ref\naszejdwa{
\ijmpa{F. Ferrari and J. Sobczyk}{11}{1996}{2213}.}--\ref\naszed{
\jgp{F. Ferrari and J. Sobczyk}{19}{1996}{287}.}. The more complicated
case of scalar fields has been treated recently
in
\ref\jmptt{F. Ferrari and J. T. Sobczyk,
{\it Jour. Math. Phys.} {\bf 41} (9) (2000), 6444,
hep-th/9909173.}.
The two and four-point amplitudes of the bosonic $\beta-\gamma$ systems
have instead been derived in
\ref\betagamma{F. Ferrari and J. T. Sobczyk,
{\it Jour. Math. Phys.} {\bf 39} (10) (1998), 5148.}.
Other useful formulas can be found in
\ref\skpss{F. Ferrari,
{\it Lett. Math. Phys.} {\bf 24} (1992), 165;
F. Ferrari,
{\it Jour. Geom. Phys.} {\bf 25} (1,2) (1998), 91;
J. Sobczyk, {\it Mod. Phys. Lett.} {\bf A8} (1993), 1153.
}.

Plane algebraic curves are described by complex coordinates which,
depending on the interpretation, can be
multivalued functions
on the complex plane {\bf C} or on the complex sphere {\bf P}$^1$.
They change their branches at the branch points
according to certain monodromy properties.
For this reason, it was thought in a first moment that
meromorphic tensors
on an algebraic curve could be easily expanded in terms of the solutions
of the related Riemann monodromy problem (see e. g.
\ref\knirev{V. G. Knizhnik,
{\it Sov. Phys. Usp.} {\bf 32}(11)
(1989) 945.}). However, this approach has encountered outstanding
obstacles in its
concrete realization. Up to now, the best way to expand meromorphic tensors
and to handle their singularity structure
is provided by the generalized Laurent series of
\naszejdwa--\naszed.
To explain their usefulness with an example, one should recall that the main
building block
in the construction of
the correlation functions of $b-c$ systems and scalar fields
is the Weierstrass kernel \weier. The latter
is a differential of the third kind which, on the curve,
plays the same role
of the Cauchy kernel on the complex plane. It is
characterized by a simple pole
in an assigned point and can be explicitly derived
with an algorithm due to Weierstrass.
However, besides the desired pole, there are also several spurious
singularities, which should be eliminated by subtracting suitable
counterterms
\zvero. A general procedure for this subtraction, which in principle
strongly depends on the form of the curve,
has been derived expanding 
the Weierstrass kernel
in generalized Laurent series.
As well, an
operator formalism on plane algebraic
curves has been established in \naszejdwa--\naszed\ .

If the plane curve is regarded as a Riemann surface,
the generalized Laurent series
can be considered as a multipoint generalization of the
Krichever-Novikov bases  \ref\kn{I. M.
Krichever and S. P. Novikov, {\it Funk. Anal. Pril.} {\bf 21} No.2 (1987), 46;
{\bf 21} No.4 (1988), 47.} much in the spirit of
\ref\raimar{R. Dick, {\it Lett. Math. Phys.} {\bf 18} (1989), 255;
\lmp{M. Schlichenmaier}{19}{1990}{151}.}.
On the other side, the modes entering in the expansion of meromorphic
tensors consist in a
finite set  of
tensors which are multivalued 
in {\bf C} (or {\bf P}$^1$). The latter have a deeper significance than
modes defined on abstract Riemann surfaces,
since they may be interpreted as correlators
of exotic conformal field theories on the complex plane (or sphere).
Theories of this kind
should be both invariant under the Virasoro group of conformal
transformations and under the monodromy group $\cal M$ of the
original curve.
The amplitudes of the class of exotic
conformal field theories corresponding to ${\cal M}=Z_n$ have been
constructed
in \ref\naszejone{F.  Ferrari, J.  Sobczyk and W.  Urbanik
{\bf 36} (1995), 3216, hep-th/9310102.}
in terms of free bosonized $b-c$ systems.
Due to the above mentioned difficulties in computing partition functions,
instead, only the $N-$point functions with $N>2$ can be derived 
if ${\cal M}$ is non abelian. The case ${\cal M}=D_n$ has been
explicitly worked out in
\ref\fercmp{\cmp{F. Ferrari}{156}{1993}{179}.}.
In this way, a concrete realization
of the more general program on holonomic quantum field theories
developed by Sato, Jimbo and Miwa in
\ref\sjm{
M. Sato, T. Miwa and M. Jimbo. Holonomic quantum fields
(Kyoto U.P. Kyoto), part I; 14 (1978) p. 223; II: 15 (1979) p. 201;
III: 15 (1979) p. 577; IV: 15 (1979) p. 871; V; 16 (1980) p.
531.} has been achieved.
The connections between generalized Laurent series and solutions of
the Riemann monodromy problem have been partly explored
in \fercmp. Some glimpses of non-commutativity in string theory
have been anticipated in 
\ref\ffbg{F. Ferrari,
{\it Int. Jour. Mod. Phys.} {\bf A9} (3) (1994), 313.}.

Despite of the above progress in understanding conformal field theories,
plane algebraic curves suffer of
some limitations.
To obtain for instance
general non-hyperelliptic Riemann surfaces of genus $g>3$, they are not
powerful enough and
it is necessary to consider algebraic curves
in the complex projective space {\bf P}$^{g-1}$.
Moreover, the smallest dimensional space in which a curve may be smoothly
embedded is ${\bf P}^3$.
In attempting to treat $b-c$ systems and scalar fields
on non-plane curves, an immediate difficulty
arises, because
the analog of the Weierstrass kernel is not known.
As explained before,
this kernel is a crucial ingredient in the derivation
of the $N-$point correlations functions and
of any other meromorphic tensor.
This problem is solved here
in the particular case of algebraic curves in {\bf P}$^3$.
A simple expression of the generalized Weierstrass kernel on these
curves is obtained and, as an application,
the correlation functions of the $b-c$ systems
and of the scalar fields are computed on a canonical curve of genus four.
Let us notice that curves in {\bf P}$^3$ include the interesting
subset containing the branched covers of Riemann surfaces.

Since it is too difficult to work on a projective variety, 
the generalization of the Weierstrass kernel has been derived in local
form after restricting oneself to
local open sets of {\bf P}$^3$,
where it is possible to use euclidean
coordinates. This is equivalent to consider affine algebraic curves
in {\bf C}$^3$.
Of course, affine algebraic curves are non-compact, but it is still possible
to relate them to compact Riemann surfaces by taking into account
also the point
at infinity and hence considering the extended complex plane
${\bf C}\cup\infty\equiv{\bf P}^1$.
A similar strategy has been followed
in the construction of
the standard Weierstrass kernel, defined in \weier\
on affine algebraic curves in {\bf C}$^2$ or on a branched cover of
${\bf P}^1$.

The material presented in this paper is divided as follows.
Section II consists in
a brief and elementary introduction to the concept of algebraic
curves. 
In Section III
the derivation of the standard Weierstrass kernel on plane curves
is reviewed and its main properties are discussed.
A formula of Weierstrass to build a differential of the second
kind on hyperelliptic curves is also presented. It is useful
to recall this formula because it
seems to have no track in the modern literature,
but it is very important whenever current-current interactions between scalar
fields have to be taken into account (see e.g. \cfft\ and the
third reference of \skpss).
Section IV contains the derivation of the generalized Weierstrass
kernel and a discussion of its main properties.
In general, we have been able to show that spurious poles can only appear
at points in which the coordinates describing the algebraic curve
in ${\bf C}^3$ become infinitely large.
In Section V
we treat the  particular case of a canonical curve of genus four in details.
As an application,
the correlation functions of the $b-c$ systems
and of the scalar fields are computed at genus four
in Section VI. These theories have been already studied using several
different
approaches,  so they
provide a good testing ground for the generalized Weierstrass kernel
of Section IV. Surface integrals over algebraic curves
are expressed as integrals in ${\bf C}^3$ in the presence of Dirac
$\delta-$functions in Appendix A. The final comments and conclusions are
presented
in Section VII.

%\ref\forsyth{A. R. Forsyth, {\it Theory of Functions of a Complex
%Variable}, Vols I and II, Dover Publications, Inc., New York, 1965.}.
%\vfill\eject
\newsec {RIEMANN SURFACES AS ALGEBRAIC VARIETIES}

Let {\bf P}$^n$ denote the complex projective space parametrized by
coordinates $\xi=\xi_0,\ldots,\xi_n$.
A homogeneous polynomial $F(\xi)$ of degree $d$ defines an algebraic
hypersurface on {\bf P}$^n$ as the locus of the points $\xi\in${\bf P}$^n$
satisfying the relation:
\eqn\basicfirst{F(\xi)=0}
The hypersurfaces are called quadrics if $d=2$, cubics if $d=3$,
quartics if $d=4$ etc.
A {\it projective algebraic variety} $\cal V$ is a subset of {\bf P}$^n$
characterized as the intersection of many hypersurfaces:
\eqn\projalgvardef{{\cal V}=\left\{\xi\in{\bf P}^n
\right|
\left. F_1(\xi)=\ldots=F_k(\xi)=0\right\}}
where $F_1,\ldots,F_k$, $k\le n$, represent a set of homogeneous polynomials
of degree $d_1,\ldots,d_k$ respectively.
If $k=n-1$, $\cal V$ describes a {\it projective algebraic curve} $C$ in
{\bf P}$^n$.
A point $p\in{\cal V}$ is {\it smooth} \ref\pgrif{P. Griffiths,
{\it Principles of Algebraic Geometry}, Wiley, New York 1978.}
if the Jacobian
matrix
\eqn\jacmat{{\cal J}=\left[{\partial(F_1,\ldots,F_k)\over
\partial(\xi_0,\ldots,\xi_n)}\right]}
has
\eqn\regcond{{\rm rank}[{\cal J}]=k}
Closely related to projective algebraic varieties are the
{\it affine algebraic varieties}:
\eqn\affalgvardef{{\cal V}_0=\left\{x\in{\bf C}^n
\right|
\left. f_1(x)=\ldots=f_k(x)=0\right\}}
$f_1,\ldots,f_k$ are polynomials
with complex coefficients and
$x=x_1,\ldots,x_n$ denotes a set of variables in the $n-$dimensional
complex space {\bf C}$^n$.
The definition of smooth points on affine varieties is analogous to that
of smooth points on projective varieties.

It is often convenient to
study projective varieties
on local open patches of {\bf P}$^n$, where it is possible to use Euclidean
coordinates. For instance, for $\sigma=0,\ldots,n$,
one may identify {\bf C}$^n$ with the following open subsets of {\bf P}$^n$
\pgrif:
\eqn\openpatchpn{U_\sigma=\left\{[\xi_0,\ldots,\xi_n]\in{\bf P}^n
\right|
\left. \xi_\sigma\ne 0\right\}}
via the homeomorphism $\phi:U_\sigma\longrightarrow{\bf C}^n$:
\eqn\homeocnpn{\phi[\xi_0,\ldots,\xi_n]=\left({\xi_0\over\xi_\sigma},
\ldots,{\xi_{\sigma-1}\over\xi_\sigma},1,{\xi_{\sigma+1}\over\xi_\sigma},
\ldots,{\xi_n\over\xi_\sigma}\right)}
with inverse
\eqn\invepcnpn{(x_1,\ldots,x_n)
\longrightarrow[x_1,\ldots,x_{\sigma},1,x_{\sigma+1},\ldots,
x_n]}
%In this way ${\cal V}_0={\cal V}\cap H_i$, where $H_i$ is the hyperplane
%orthogonal with respect to $\xi=0$.
To simplify the notations, we will consider hereafter only the case
$\sigma=0$. In doing that, there is no loss of generality, since none
of the variables $\xi_1,\ldots,\xi_n$ plays a privileged role in the present
context.
Clearly, upon the identifications
\eqn\idents{x_1={\xi_1\over\xi_0},\ldots,x_n={\xi_n\over\xi_0}}
and
\eqn\identt{f_a(x_1,\ldots,x_n)\equiv F_a(1,x_1,\ldots,x_n)\qquad\qquad
a=1,\ldots,k}
the restriction of a projective variety ${\cal V}$ on $U_0$ is equivalent to
an affine variety on {\bf C}$^n$ associated to the system of
equations:
\eqn\assalgequ{f_a(x_1,\ldots,x_n)=0\qquad\qquad a=1,\ldots,k}

In the special case $k=n-1$, a projective curve and the affine curve
related to it via the homeomorphism \homeocnpn\ differ only
by a finite number of points ``at infinity''. The latter are
determined by the conditions $F_1=\ldots=F_{n-1}=\xi_0=0$, in which
$\phi$ is no longer defined.
To show that, we consider a projective algebraic curve
\eqn\projalgcurdef{C=\left\{\xi\in{\bf P}^n
\right|
\left. F_1(\xi)=\ldots=F_{n-1}(\xi)=0\right\}}
Since the $F_a$ are homogeneous
polynomials of degree $d_a$, $a=1,\ldots,n-1$, they can be written as follows:
\eqn\homogpoly{F_a(\xi)=\sum_{i_1,\ldots,i_n=0}^{d_a}
(\xi_1)^{i_1}\ldots(\xi_n)^{i_n}(\xi_0)^{d_a-(i_n+\ldots+i_1)}}
If $\xi_0=0$, the system of algebraic equations defining $C$ is trivial
unless the condition $d_a-(i_n+\ldots+i_1)=0$ is fulfilled. Thus,
the points on the curve corresponding to $\xi_0=0$ are given by the residual
system of equations:
\eqn\ressys{
\sum_{i_1,\ldots,i_{n-1}=0}^{d_a}
(\xi_1)^{i_1}\ldots(\xi_{n-1})^{i_{n-1}}(\xi_n)^{d_a-(i_{n-1}+\ldots+i_1)}=0
}
for $a=1,\ldots,n-1$.
The above relations describe the intersection of $n-1$ hypersurfaces in
{\bf P}$^{n-1}$, which, by B\'ezout theorem, contains a finite number of
$d_1d_2\cdots d_{n-1}$ common points as desired.
\smallskip
Algebraic curves are particularly important in the study of compact
{\it Riemann surfaces}. The simplest representation of a Riemann surface
$\Sigma_g$ of genus $g$ is in terms of
{\it plane projective algebraic curves}, or hypersurfaces in
{\bf P}$^2$
associated to a single algebraic equation of the kind:
\eqn\planalgcurdef{F(\xi_0,\xi_1,\xi_2)=0}
A plane curve is said {\it non-singular} or {\it regular}
provided the condition below is never
verified:
\eqn\placurnonsing{F={\partial F\over\partial \xi_0}=
{\partial F\over\partial \xi_1}=
{\partial F\over\partial \xi_2}=0}
Modulo conformal transformations, any compact Riemann surface coincides with
a plane projective algebraic curve.\smallskip
A plane curve $C$ can be projected
from a point $p\not\in C$ 
onto a complex line $L$ in {\bf P}$^2$,
which does not contain $p$. After a linear change of coordinates
one may take $p=[0,0,1]$
and $L=\{\xi\in{\bf P}^2|\xi_2=0\}$. The result of the projection
is a representation of the Riemann surface as
a {\it branched cover} (or $n-${\it sheeted covering})
of {\bf P}$^1$, whose points are given by the
zeros of the complex polynomial $f(x_1,x_2)=F(1,x_1,x_2)$ in the Euclidean
coordinates $x_1=\xi_1/\xi_0$ and $x_2=\xi_2/\xi_0$.
Solving the equation
\eqn\bcovering{f(x_1,x_2)=0}
with respect to $x_2$, one obtains a multivalued function $x_2(x_1)$
of $x_1\in{\bf P}^1$. The finite {\it branch points} of $x_2(x_1)$ satisfy the
relations:
\eqn\brapoidef{f(x_1,x_2)=f^{(0,1)}(x_1,x_2)=0}
where we have used the notation
\eqn\pardernot{f^{(m,n)}(x_1,x_2)\equiv {\partial^m\over\partial x_1^m}
{\partial^n\over\partial x_2^n} f(x_1,x_2)}
Starting from a multivalued function, a Riemann surface $\Sigma_g$ can be
constructed in terms of sheets and cross-cuts.
In analogy with eq. \placurnonsing, a cover of {\bf P}$^1$
is non-singular if the
following identities are never satisfied simultaneously:
\eqn\bracovnonsing{f={\partial f\over\partial x_1}=
{\partial f\over\partial x_2}=0}
\smallskip

The double-sheeted coverings of {\bf P}$^1$, or {\it
hyperelliptic curves}, have very special properties with respect to the other
curves. 
Their polynomial equation can always be reduced
to the following one:
\eqn\hypcurdef{x_2^2=P_{2g+2}(x_1)}
where $P_{2g+2}(x_1)$ is a polynomial of degree
$2g+2$ in $x_1$ with complex coefficients.
Exploiting the group $SL(2,{\bf C})$ of automorphisms
of the sphere {\bf P}$^1$, the number of independent coefficients
reduces to $2g-1$.
Since a general Riemann surface of genus $g>2$ depends
on $3g-3$ complex parameters, called the {\it moduli}, not all Riemann
surfaces can be hyperelliptic.
As a matter of fact, hyperelliptic curves
form only a subset of dimension $2g-1$ in the moduli space \fakra.
On the other side,
plane algebraic curves of genus $g\ge 3$ may be also non-hyperelliptic,
but the number of independent parameters which is possible to accommodate
in the defining polynomial $F(\xi)$ of eq.
\planalgcurdef\ is less than $3g-3$ when $g>3$.
For this reason, in order to construct general non-hyperelliptic Riemann
surfaces, one
usually considers {\it canonical curves} embedded in {\bf P}$^{g-1}$.

A canonical map $\varphi$ from a Riemann surface $\Sigma_g$
to the projective space ${\bf P}^{g-1}$ can be established with the help of
a basis of holomorphic differentials $\omega_1,\ldots,\omega_g$
on $\Sigma_g$ as follows:
\eqn\canmap{\varphi:p\in\Sigma_g\longrightarrow[\omega_1,\ldots,\omega_g]\in
{\bf P}^{g-1}}
If $\Sigma_g$ is a general non-hyperelliptic Riemann surface,
$C=\varphi(\Sigma_g)$ is called a canonical curve.
For instance, a canonical curve of genus $g=3$ is a plane projective algebraic
curve of degree four. The canonical curve with $g=4$ is instead given by the
complete intersection of a quadric and a cubic in 
${\bf P}^3$. For $g=5$, the canonical curve is a complete intersection
of three quadrics in ${\bf P}^4$ etc.
A more detailed classifications of canonical curves
together with a thorough 
discussion of some exceptional cases can be found in \cenr.
\smallskip
Here we will limit ourselves to algebraic curves in {\bf P}$^3$, which
is the smallest dimensional
projective space in which a curve can be smoothly embedded.
In fact, smooth embeddings are not possible in {\bf P}$^2$,
so that
plane algebraic curves are affected by singularities
at isolated points
\foot{However, a curve $C$ in {\bf P}$^3$ may always be projected in
{\bf P}$^2$ in
such a way that the
resulting plane algebraic curve has only ordinary double points.}.
Let $F(\xi)$ and $G(\xi)$,
$\xi=\xi_0,\xi_1,\xi_2,\xi_3$ be two homogeneous polynomials of
degrees $d_F$ and $d_G$
respectively, whose zeros define two hypersurfaces in {\bf P}$^3$.
We are mainly interested in situations in which the intersection of
the two hypersurfaces is {\it complete}, i.e. they
meet in a single curve $C$, whose points are given
by the following system of algebraic equations:
\eqn\prosysalgeqs{F(\xi_0,\xi_1,\xi_2,\xi_3)=G(\xi_0,\xi_1,\xi_2,\xi_3)=0}
We  also assume that all the points of $C$ are smooth.
The genus $g$ of $C$ is then given by:
\eqn\comintgen{2g-2=d_Fd_G(d_F+d_G-4)}

Supposing that $\xi_0\ne 0$ and using the mapping
\homeocnpn, we obtain the affine algebraic variety in {\bf C}$^3$
associated to the polynomial equations:
\eqn\sysalgeqs{f(x_1,x_2,x_3)=g(x_1,x_2,x_3)=0}
Apart from the exceptional cases described in
\cenr, the above relations describe
the complete intersections of two hypersurfaces in {\bf C}$^3$.
One may also view $f(x)$ and $g(x)$ as polynomials
of degrees $m\le d_F$ and $n\le d_G$ respectively in the
variable $x_3$. Eliminating the latter, the 
resultant $R(x_1,x_2)$
is a polynomial in $x_1$ and $x_2$ of degree $mn$ and
the equation
\eqn\respro{R(x_1,x_2)=0}
represents the projection of the curve \sysalgeqs\ on {\bf C}$^2$.
Let us notice that, with respect to the case of plane curve, this
projection is not unique. For instance, it is possible to
eliminate from \sysalgeqs\ the variable $x_2$ instead of $x_3$.
In the latter case one obtains a different resultant $R'(x_1,x_3)$
and a different projection onto {\bf C}$^2$:
\eqn\resprobis{R'(x_1,x_3)=0}

To \sysalgeqs\ one can also associate a compact Riemann
surface $\Sigma_g$ of genus $g$ constructed as a ramified covering of
{\bf P}$^1$.
For example,
in the neighborhood of a point where the Jacobian
\eqn\jactt{J^1(x)={\partial f(x)\over \partial x_2}
{\partial g(x)\over \partial x_3}-
{\partial g(x)\over \partial x_2}{\partial f(x)\over \partial x_3}}
is different from zero,
$x_1\in{\bf P}^1$ becomes a good local coordinate and it is possible
to solve the system of algebraic equations
\sysalgeqs\ with respect to $x_2$ and $x_3$.
In this way, one obtains two multivalued functions
$x_2(x_1)$ and $x_3(x_1)$, whose analytic continuation on the complex
line {\bf P}$^1$ defines a Riemann surface.
The ramification points are those in which the condition
$J^1(x)=0$ is satisfied.
We stress again the difference with respect to
plane curves, because now there 
are two possible representations of $\Sigma_g$ in terms of
branched covers of ${\bf P}^1$, corresponding to
eqs. \respro\ and \resprobis. 

Finally, we show that eqs. \sysalgeqs\ includes the branched covers
of Riemann surface as a particular sub-case.
Indeed, let us consider a compact Riemann surface represented
as a branched cover of ${\bf P}^1$
associated to the vanishing of a polynomial:
\eqn\bracovcs{f(x_1,x_2)=0}
Any Riemann surface of genus
$g$ can be
mapped into  a branched cover of this kind
in a limited region of its moduli space.
Solving eq. \bracovcs\ with respect to $x_2$,
the Riemann surface $\Sigma_g$ is parametrized as a curve
${\bf C}^2$ with coordinates $(x_1,x_2(x_1))$.
Starting from such coordinates it is still possible to realize a branched
cover $\tilde\Sigma$ of $\Sigma_g$ by requiring that:
\eqn\bracovrs{g(x_2,x_3)=0}

%This plane curve can be further projected on a complex line $L$
%as explained above. Thus algebraic curves in {\bf P}$^3$
%may be associated to a branched cover
%of {\bf P}$^1$ like plane curves.
%Let us choose for example
%$L=\{\xi\in{\bf P}^3|\xi_2=\xi_3=0\}$. The
%%system of equations \sysalgeqs\
%with respect to the variables $x_2$ and $x_3$ Explicitly, this is done
%by computing its resultant,  which
%%will become multivalued functions
%of $x_1\in{\bf P}^1$. Depending on the way in which the projection is made,
%one has four different possibilities:
%\eqn\posone{x_3=x_3^f(x_2(x_1),x_1),\qquad x_2=x_2^g(x_1)}
%\eqn\postwo{x_3=x_3^g(x_2(x_1),x_1),\qquad x_2=x_2^f(x_1)}
%\eqn\posthr{x_2=x_2^f(x_3(x_1),x_1),\qquad x_3=x_3^g(x_1)}%
%%\eqn\posfou{x_2=x_2^g(x_3(x_1),x_1),\qquad x_3=x_3^f(x_1)}
%where $x_i^f$, $i=2,3$, means that $x_i$ has been derived solving the
%algebraic equation $f=0$. An analogous definition is valid for
%$x_i^g$.

%The theory of conformal free fields on n-sheeted coverings of ${\bf P}^1$
%has been already treated in details. In the next Section the problem
%will be studied on Riemann surfaces embedded in ${\bf P}^3$ of the kind:
%\eqn\algembpt{C=\left\{\xi\in {\bf P}^3\right|\left.
%F(\xi)=G(\xi)=0\right\}}
\newsec{THE WEIERSTRASS KERNEL}
The Cauchy kernel
\eqn\cauker{K_C(x_1,x_2)={dx_1\over x_1-x_2}}
plays a fundamental role in the construction of the amplitudes of conformal
field theories on the complex plane {\bf C}.
This kernel has the following two properties \zvero:
\item{1)} As a function of $x_1$, the kernel is a meromorphic differential
with two simple poles in $x_1=x_2$ and in $x_1=\infty$. The residues are
$+1$ and $-1$ respectively.\medskip
\item{2.)} As a function of $x_2$, the kernel is a meromorphic function
 with a simple pole in $x_2=x_1$ and a single zero in $x_2=\infty$.
\medskip
A meromorphic function with a single pole cannot exist on
a compact two dimensional manifold. For this reason, usually it is only
required that an analogue $K(p,q)dp$ of the Cauchy kernel on a Riemann surface
$\Sigma_g$
should have the following asymptotic behavior in a neighborhood of the point
$q$\zvero:
\eqn\propess{K(p,q)\sim{dp\over p-q}+A(p,q)}
where $A(p,q)$ is finite at $p=q$.\smallskip
A kernel with the above property can be constructed on $n-$sheeted coverings
of {\bf P}$^1$ using a well-known algorithm of Weierstrass \weier.
To this purpose, let us consider a Riemann surface $\Sigma_g$ represented
as the locus of points defined by eq. \bcovering. Solving this equation
for $x_2$, one
obtains a multivalued function $x_2(x_1)$, while $x_1$ is a $d-$degree mapping
$x_1:\Sigma_g\longrightarrow{\bf P}_1$
from the Riemann surface to the projective sphere {\bf P}$^1$.
Any point $p\in \Sigma_g$ is in a $1-1$ correspondence
with a point on the branched cover
$x(p)\equiv x_1(p),x_2(x_1(p))$.
\smallskip
Let us put now
$y(q)\equiv x_1(q),x_2(x_1(q))$, where $q\ne p$.
Then an analogue of the Cauchy kernel on the branched cover is the following
Weierstrass kernel:
\eqn\weiker{K_W(x,y)={f(y_1,x_2)\over(x_2-y_2)f^{(0,1)}(x_1,x_2)}
\enskip{dx_1\over x_1-y_1}}
To study its behavior in the limit $x_1\longrightarrow y_1$,
one can expand $f(y_1,x_2)$
in series of Taylor near the point $x_2=y_2$:
\eqn\ytexp{f(y_1,x_2)\sim f(y_1,y_2)+f^{(0,1)}(y_1,y_2)(x_2-y_2)+\ldots}
Since $f(y_1,y_2)=0$, it is clear that the Weierstrass kernel 
satisfies
requirement \propess:
\eqn\weiksimple{K_W(x,y)\sim{dx_1\over x_1-y_1}}
The next term in the Taylor expansion gives in fact a contribution
proportional to
\eqn\finterm{A(x_1,y_1)={f^{(0,2)}(y_1,y_2)f^{(1,0)}(x_1,x_2)\over
f^{(0,1)}(x_1,x_2)}dx_1}
which is finite when $x_1=y_1$.

One should also consider those
points $p\ne q\in \Sigma_g$ which, on the
branched cover, correspond to coordinates $x(p)$
and $y(q)$ such that
$x_1(p)=x_1(q)$ and $x_2(x_1(p))\ne x_2(x_1(q))$. 
In this case, there is no spurious pole in the Weierstrass kernel despite the
presence of the factor $(x_1-y_1)^{-1}$.
As a matter of fact, expanding the function
$f(y_1,x_2)$ in series of Taylor
in $y_1$:
\eqn\tayzeta{f(y_1,x_2)=f(x_1,x_2)+f^{(1,0)}(x_1,x_2)(x_1-y_1)+\ldots}
and substituting in \weiker\ one obtains
\eqn\anompoint{K_W(x,y)\sim
{f^{(1,0)}(x_1,x_2)\over (x_2-y_2)f^{(0,1)}(x_1,x_2)}dx_1}
which is finite since
$x_2-y_2\ne 0$ by hypothesis.

One may also check that
the Weierstrass kernel has no spurious poles at the
branch points if the curve is regular.
To this purpose it is possible to exploit the relation
\eqn\bccift{{dx_1\over f^{(0,1)}(x_1,x_2)}={dx_2\over f^{(1,0)}(x_1,x_2)}}
which is a consequence of the implicit function theorem \pgrif.
Let us now consider a
branch point $x_1=x_1(a)$ of {\it multiplicity}
$\nu$, where $\nu$ is an integer and $a\in \Sigma_g$.
We suppose for the moment that $x_2(x_1(a))$ is finite.
Near the branch point we choose a good local coordinate $t$ such that:
\eqn\loccoor{x_1-x_1(a)=t^\nu}
Since $x_2(x_1(a))$ is not divergent, its approximate expansion in powers of
$t$ will look as follows:
\eqn\xtloccoord{x_2\sim\alpha_0+\alpha_1 t+\alpha_2 t^2+
\ldots}
with $\alpha_{0,1,2}$ being constants.
Thus $dx_1\sim \nu t^{\nu-1}dt$ and $dx_2\sim dt$.
Remembering that
the function $f^{(1,0)}(x_1,x_2)$
does not vanish at a branch point
due to the regularity hypothesis \regcond, we find
that near $x_1(a)$ eq. \bccift\ is approximated by:
\eqn\absspupolbp{{dx_1\over f^{(0,1)}(x_1,x_2)}\sim dt}
This shows that the zeros of $f^{(0,1)}(x_1,x_2)$
are absorbed by the corresponding zeros of the
differential $dx_1$, so that the Weierstrass kernel
cannot be singular at the finite branch
points. If $x_2$ has a pole of order $k$ near
a branch point, instead, it is always possible to perform
the change of variables
\eqn\chavar{x_2'=(x_1-x_1(a))^k x_2}
In the new coordinates, $x_2'$ remains finite at the branch point
$x_1(a)$ and the above demonstration applies again.

Of course, on a Riemann surface a single pole is not allowed, so that
the Weierstrass kernel \weiker\ must contain also other spurious
poles both as a function of $x_1$ and $y_1$.  Typically, they appear
whenever the variables $x_1,x_2$ and $y_1,y_2$ become
infinitely large.
A general procedure to subtract these spurious poles while keeping the
property \propess\ has been already discussed in ref.  \naszejdwa--\naszed\.
The amplitudes of some free conformal
field theories on branched covers of {\bf P}$^1$ have been  constructed
in.

To conclude this Section, we present a beautiful formula of Weierstrass
\weier\
to construct a second kind differential on an hyperelliptic curve.
Thus, we consider curves of the kind \hypcurdef,
where
\eqn\hyppol{P_{2g+2}(x_1)=A_0+A_1x_1+\ldots+A_{2g+2}x^{2g+2}}
Let us define the function:
$$
R(x_1,x'_1)=A_0+{1\over 2}A_1(x_1+x_1')+A_2x_1x_1'+
{1\over 2}A_3(x_1+x_1')x_1x_1'+A_4(x_1x_1')^2+
$$
\eqn\rxxp{+{1\over 2}A_5(x_1+x_1')(x_1x_1')^2+\ldots
+{1\over 2}A_{2g+1}(x_1+x_1')(x_1x_1')^g+A_{2g+2}(x_1x_1')^{g+1}}
Clearly, $R(x_1,x'_1)$ satisfies the property
\eqn\propppp{R(x_1,x_1')=R(x_1',x_1)}
Moreover, if $x_1=x_1'$, one has that:
\eqn\prere{R(x_1,x_1')=P_{2g+2}(x_1)\qquad\qquad
\left. {\partial R(x_1,x_1')\over\partial x_1}\right|_{x_1=x_1'}=
{1\over 2}{\partial
P_{2g+2}(x_1)\over \partial x_1}}
The differential of the second kind is given by:
\eqn\seckindif{\tau_{x_1'}(x_1)=
-{y(x_1)y'(x_1)+R(x_1,x'_1)\over2(x_1-x_1')^2y(x_1)y'(x_1)}dx_1}
where $y'(x_1)={\partial y(x_1)\over\partial x_1}$.
It is easy to show using the properties \propppp-\prere\ that
$\tau_{x_1'}(x_1)$ has only a pole of the second order in
$x_1=x_1'$ as it should be.
Unfortunately it is not simple to extend the elegant formula \seckindif\
to the $Z_n$ symmetric curves, not to mention the general plane algebraic
curves of eq. \bcovering\ or the even more complicated curves \sysalgeqs.
\newsec{GENERALIZED WEIERSTRASS KERNELS}
In this Section we construct analogues of the Weierstrass kernel on
 affine algebraic curves defined
by the system of equations \sysalgeqs.
Even if it will not be strictly necessary, we suppose to fix the ideas
that the 
intersection of the two hypersurfaces in \sysalgeqs\ is complete and
gives as a result an algebraic curve $C$ which coincides, modulo conformal
transformations, with a Riemann surface $\Sigma$.

The following two different cases can formally be treated in the same way.
On one side, $x_1,x_2,x_3$ may be interpreted as coordinates
in {\bf C}$^3$, so that $\Sigma$ is not compact due
to the absence of the points at infinity. Alternatively, the vanishing
of the polynomials \sysalgeqs\ can be associated to a ramified covering
of {\bf P}$^1$ as we have seen in Section II and $\Sigma$ is
a compact Riemann surface of
genus $g$.

Let us consider as in the previous Section
two different points $p,q\in \Sigma$.
On the algebraic curve $C$ they correspond to coordinates
$x(p)=x_1(p),x_2(p),x_3(p)$ and $y(q)=y_1(q),y_2(q),y_3(q)$.
A possible analogue of the Weierstrass kernel on $C$ is given by:
\eqn\ksymm{K_{\rm sym}(x,y)={1\over 3}\sum_{i=1}^3
{N^i(x,y)\over J^i(x)}
{dx_i\over (x_1-y_1)(x_2-y_2)(x_3-y_3)}}
where
\eqn\jix{J^i(x)=\epsilon^{ikl}{\partial f(x)\over\partial x_k}
{\partial g(x)\over\partial x_l}}
and
\eqn\nixyone{N^1(x,y)=f(y_1,y_2,x_3)g(x_1,y_2,x_3)-
f(x_1,y_2,x_3)g(y_1,y_2,x_3)}
\eqn\nixytwo{N^2(x,y)=f(x_1,y_2,y_3)g(x_1,x_2,y_3)-
f(x_1,x_2,y_3)g(x_1,y_2,y_3)}
\eqn\nixythree{N^3(x,y)=f(y_1,x_2,y_3)g(y_1,x_2,x_3)-
f(y_1,x_2,x_3)g(y_1,x_2,y_3)}
Here $\epsilon^{ikl}$ denotes the completely antisymmetric tensor in
three dimensions with the convention $\epsilon^{123}=1$.
We note that the variables $x,y$ and the functions $f$ and $g$ enter
symmetrically in the expression of the kernel \ksymm, as it should be
since none of them plays a privileged role in the definition of the
algebraic curve.
The symmetry under the exchange of $f$ and $g$ is
also related to the freedom of projecting the curve in the two
possible
ways shown by eqs. \respro\ and \resprobis.

Equivalent kernels can be obtained starting from $K_{\rm sym}(x,y)$ and adding
differentials in such a way that the behavior near the singularity in
$x=y$ remains unchanged. For instance, exploiting the identities:
\eqn\basrel{{dx_1\over J^1}={dx_2\over J^2}=
{dx_3\over J^3}}
which are a consequence of the implicit function theorem
and using the fact that the numerators $N^1,N^2,N^3$ differ each other
by functions that vanish in $x=y$,
one may derive the following kernel:
\eqn\fker{K(x,y)={N^1(x,y)\over J^1(x)}
{dx_1\over (x_1-y_1)(x_2-y_2)(x_3-y_3)}}
The above kernel is less symmetric than $K_{\rm sym}(x,y)$, but has a more
compact expression. 
In the particular case
of a branched cover of a plane curve,
in which the algebraic equations
\sysalgeqs\ assume the form \bracovcs\ and \bracovrs,
$K(x,y)$ is simply given by:
\eqn\simwgwk{
K(x,y)=-{f(x_1,y_2)g(y_2,x_3)\over G_{x_3}(x)F_{x_2}(x)}
{1\over (x_1-y_1)(x_2-y_2)(x_3-y_3)}
}

Let us investigate the behavior
of $K(x,y)$ near $x=y$.
To this purpose, it is sufficient to expand $N^1(x,y)$ with respect to
$y$ at the point $x$. At the leading order:
\eqn\noexpyx{N^1(x,y)\sim J^1(x)\Delta x_2\Delta x_3}
where $\Delta x_i=y_i-x_i$ is very small by hypothesis.
We note that, in principle, there are also
contributions proportional to $\Delta x_1$
in the expansion of $N^1(x,y)$. In order to obtain eq.
\noexpyx, $\Delta x_1$
has been expressed in terms of $\Delta x_2$ and $\Delta x_2$
with the help of \basrel.
Substituting \noexpyx\ in \fker\ we find that $K(x,y)$ satisfies
the property \propess\ as desired:
\eqn\fkerasym{K(x,y)\sim{dx_1\over (x_1-y_1)}}
A similar calculation for $K_{\rm sym}(x,y)$ gives the result
\eqn\ksymmasym{K_{\rm sym}(x,y)\sim{1\over 3}
\sum_{i=1}^3{dx_i\over (x_i-y_i)}}
In this case one should remember that
not all variables $x_1,x_2,x_3$ are independent due to
the relations \sysalgeqs, so that eq.~\ksymmasym\ must be worked out further.
Assuming for instance that $x_1$ is a good local coordinate in a
neighborhood
of the point $p$, one can solve the system of algebraic equations
\sysalgeqs\ with respect to the remaining variables, so that
$x_j=x_j(x_1)$ for $j=2,3$.
In the same way $y_j=x_j(y_1)$ for $j=2,3$ in a neighborhood of $q$, with
$y_1=x_1+\Delta x_1$.
As a consequence
\eqn\varie{dx_j(x_1)={dx_j\over dx_1}dx_1\qquad\qquad
x_j(y_1)-x_j(x_1)\sim -{dx_j\over dx_1}\Delta x_1\qquad\qquad j=2,3}
Using the above relations in \ksymmasym, one obtains:
\eqn\finksymmasym{K_{\rm sym}(x,y)\sim{dx_1\over (x_1-y_1)}}
Thus, also the kernel \ksymm\ has the requested behavior
near the pole in $x=y$.
\smallskip
Besides the required simple pole in $x=y$, the
kernels \ksymm\ and \fker\ have also spurious poles, which have to be
controlled and suitably
subtracted in order to construct physical correlation functions
with desired singularities.
The study of these spurious poles will be the subject of the rest
of this Section.

Since the structure of the kernels
\ksymm\ and \fker\ consists in ratios of polynomials
of  $x$ and $y$, their possible divergences may only occur
at the zeros of the denominators
$J^i(x)(x_1-y_1)(x_2-y_2)(x_3-y_3)$ or at the infinities of the
numerators
$N^i(x,y)$, $i=1,2,3$.

First of all, we consider the zeros of $J^1(x)$.
The cases in which $i=2,3$ can be treated
in an analogous way.
Let $x(a)=(x_1(a),x_2(a),x_3(a))$, $a\in \Sigma$, be
a point on $C$ for which $J^1(x)=0$.
In $x(a)$ the system of algebraic equations
\sysalgeqs\ becomes no longer invertible with respect to
$x_2$ and $x_3$. Given a good local coordinate $t$ in a neighborhood
of $x(a)$, this implies that
\eqn\ttt{x_1-x_1(a)=t^\lambda}
and
\eqn\tttt{x_2=\alpha_0+\alpha_1 t^\mu+\ldots\qquad\qquad
x_3=\beta_0+\beta_1 t^\nu+\ldots}
where $\lambda$ is an integer containing the integers $\mu$ and $\nu$
as sub-factors. In eq. \tttt\ we suppose
that $x_2$ and $x_3$ do not diverge in $a$, so that
$\lambda,\mu,\nu$ are all positive. If not, it is always possible to
perform in eq. \sysalgeqs\ a change of variables
$x_2,x_3\rightarrow x_2',x_3'$ similar to that of eq. \chavar\
to make the new variables $x_2',x_3'$ finite in $a$.
We note that the monodromy properties
of $x_2$ and $x_3$ may be in general different, so that $\mu$ and $\nu$
need not to be equal.
Near $x(a)$ the relations \basrel\ become:
\eqn\brxatp{{\lambda t^{\lambda-1}dt\over J^1(x)}=
{\mu t^{\mu-1}dt\over J^2(x)}=
{\nu t^{\nu-1}dt\over J^3(x)}}
where $x$ is a function of $t$ given by eqs. \ttt\ and \tttt.
At this point it is possible to invoke the regularity condition
\regcond, which assures that $J^2(x),J^3(x)\ne 0$ in $x=x(a)$.
Since $\mu,\nu\le \lambda$ it is easy to see from \brxatp\ that
the zeros of $J^1(x)$ are absorbed by the corresponding zeros of the
differential $dx_1$. As a consequence, there are no spurious divergences
at these points.

One can also verify that there are no poles when
$(x_1-y_1)(x_2-y_2)(x_3-y_3)=0$ apart from the
one in $x(p)=y(q)$, which is related
to the required singularity in $p=q$.
Spurious poles of this kind may in principle occur if
two different points $p,q\in\Sigma$ correspond on the algebraic
curve to coordinates
$x(p)$ and $x(q)$ characterized by the fact that some of their components,
but not all, coincide (for instance $x_1=y_1,x_2=y_2$ and $x_3\ne y_3$).
The proof that no spurious divergence arises in this case is straightforward
and will not be reported here.

In conclusion, the kernels \ksymm\ and \fker\
diverge only at the poles of the numerators $N^i(x,y)$.
In general, the latter are located at the points in which the variables
$x$ and $y$ become very large.

\newsec{THE CASE OF GENERAL NON-HYPERELLIPTIC CURVES OF GENUS FOUR}
A general non-hyperelliptic algebraic curve of genus four is given by
the complete intersection of a quadric with a cubic as mentioned in
Section II. In fact,
putting $d_F=3$ and $d_G=2$ in \comintgen, one obtains exactly $g=4$.
For instance. one can choose in eqs. \prosysalgeqs\
$G(\xi)$ as follows:
\eqn\pggf{G(\xi_1,\xi_2,\xi_3,\xi_4)=\xi_0\xi_1-\xi_2\xi_3}
while $F(\xi)$ is a homogeneous polynomial of degree three.
In affine coordinates $\xi_i=\xi_i/\xi_0$, $i=1,2,3$, $F(\xi)$ and
$G(\xi)$ are replaced by:
\eqn\ggf{G(1,x_1,x_2,x_3)=g(x_1,x_2,x_3)=x_2x_3-x_1}
\eqn\fgf{F(1,x_1,x_2,x_3)=f(x_1,x_2,x_3)=
x_3^3+h_1(x_1,x_2)x^2_3+h_2(x_1,x_2)x_3+h_3(x_1,x_2)
}
where
\eqn\hixx{h_i(x_1,x_2)=\sum_{k,l=0\atop k+l\le 3}^ia_{kl}^{(i)}
x_1^kx_2^l\qquad\qquad
i=1,2,3}
and the $a_{kl}^{(i)}$ are complex coefficients.
We note that the polynomial $f(\xi_1,\xi_2,\xi_3)$ has been
ordered according to the different powers of $x_3$. This is just a convention
which does not reflect any special role of $x_3$. In the same
way one could order $f(\xi_1,\xi_2,\xi_3)$ with respect to the
powers of $x_1$ or $x_2$.
All necessary ingredients to construct the Weierstrass kernels
\ksymm\ and \fker\ are derived in a straightforward way
substituting eqs. \ggf\ and \fgf\ in \jix-\nixythree.

In the following we will assume that $x_1$ is a good local coordinate of the
curve. The cases in which $x_1\in${\bf C} or $x_1\in${\bf P}$_1$ can be
formally treated in the same way.
Accordingly, we solve the system of algebraic equations
\eqn\gfsaeq{f(x_1,x_2,x_3)=g(x_1,x_2,x_3)=0}
with respect to $x_1$\foot{There is no loss of generality in doing that. If
one wishes to study the algebraic curve in the neighborhood of a branch point,
where $x_1$ is no longer a good coordinate, it is always possible to
perform a conformal transformation and to consider $x_1$ as a function of
$x_2$ or $x_3$.}. As a result,
one obtains two multivalued functions
$x_2(x_1)$ and $x_3(x_1)$. Due to eq. \ggf, both $x_2(x_1)$ and $x_3(x_1)$
share the same monodromy properties. 
They define a Riemann surface
$\Sigma_4$ constructed in terms of sheets.
It is easy to check that
$x_2(x_1)$ and $x_3(x_1)$ have six branches, so that $\Sigma_4$ consists
of six sheets glued together at the branch lines.
The computation of the resultant $R(x_1,x_2)$ of the two algebraic
equations \gfsaeq\ is straightforward and gives the following
plane curve equation associated to $\Sigma_4$:
\eqn\gfplacureqn{x_1^3+h_1(x_1,x_2)x_1^2x_2+h_2(x_1,x_2)x_1x_2^2+
h_3(x_1,x_2)x_2^3=0}
which indeed describes a curve of genus four.

Due to the peculiar role played by the variable $x_1$,
it is natural to consider the kernel $K(x,y)$ instead of the more
symmetric one of eq. \ksymm. Inserting eqs. \ggf\ and \fgf\ in
\fker\ one obtains:
\eqn\gfker{K(x,y)=-{dx_1\over(x_1-y_1)J^1(x)}\left[
{x_3f(y_1,y_2,x_3)\over (x_3-y_3)}+{y_2f(x_1,y_2,x_3)\over
(x_2-y_2)}\right]}
where we should remember that we are dealing with multivalued
functions
$x_j=x_j(x_1)$ and $y_j=y_j(y_1)$ for $j=2,3$.
With respect to the general formula \fker, in this case
some simplifications
have been possible
because of the particular form of
$g(x_1,x_2,x_3)$.
One can check that the above kernel has the desired pole when
$x_i=y_i$ for $i=1,2,3$. If one wishes to study the kernel \gfker\ in the
neighborhood of a branch point where $J^1(x)=0$, it is possible to perform
the conformal transformation $x_1=x_1(x_2)$.

Let us now concentrate
on the spurious divergences of $K(x,y)$.
From the previous Section, we know that they
may only occur at the infinities of the variables
$x_i,y_i$, $i=1,2,3$. In this case, the situation is made simpler
by the fact that $x_2(x_1)$ and
$x_3(x_1)$ have no poles
for finite values
of $x_1$. To show that, let us imagine that a point $a\in \Sigma_4$
corresponds on the algebraic curve to a point $x_1(a)$ where
$x_3$ has a pole of order $s$:
\eqn\finsin{x_3(x_1)\sim (x_1-x_1(a))^{-s}}
Due to eq. \ggf, the function $x_2(x_1)$ has a zero of the same
order at the same point:
\eqn\finzer{x_2(x_1)\sim (x_1-x_1(a))^{s}}
Thus, in the limit $x_1=x_1(a)$ the polynomials
$h_i(x_1,x_2)$, $i=1,2,3$, may be replaced 
by suitable constants
$A_1,A_2,A_3$ neglecting
higher order terms
in $x_1-x_1(a)$.
Hence, eq. \fgf\ is approximated by:
\eqn\appfgf{x_3^3+A_1x^2_3+A_2x_3+A_3=0}
Clearly, the above equation has no solutions if $x_3=\infty$.
Analogously, since $x_3$ and $x_2$ enter in eq. \gfsaeq\ symmetrically,
it is possible to verify that
$x_2(x_1)$ has no divergences for finite values of $x_1$.

To study the singularities of $x_2(x_1)$ and $x_3(x_1)$ at infinity,
it is convenient to introduce the new variable
$x_1'=x_1^{-1}$.
Let us now suppose that $x_2$ and $x_3$ 
have the following behavior near $x_1'=0$:
\eqn\posinf{x_2=\alpha (x_1')^{s}+\ldots\qquad\qquad x_3=\beta
(x_1')^{-1-s}+\ldots}
The second
relation \posinf\ is again a consequence of \ggf.
Substituting the ansatz \posinf\ in \fgf, it is easy to verify that
the latter equation is satisfied
only if $s=-1$ or $s=0$.
In the first case, there are
three
branches of $x_2$ and $x_3$ such that $x_2$ diverges:
\eqn\thrbraone{x_2\sim {1\over x_1'}\qquad\qquad x_3\sim{\rm const}}
If $s=0$, instead, there are other three branches in which
$x_3$ becomes singular:
\eqn\thrbratwo{x_2\sim {\rm const}
\qquad\qquad x_3\sim{1\over x_1'}}

At this point, we are ready to discuss the spurious poles of the kernel
\gfker.
As a meromorphic differential in $x_1$, $K(x,y)$ has three simple poles
in $x_1=0$.
The latter occur in the three branches of $x_2$ and $x_3$ where eq.
\thrbratwo\ is satisfied.
In each of these branches, $K(x,y)$ has residue $-{1\over 3}$:
\eqn\gfresinf{K(x,y)=-{1\over 3}{dx_1'\over x_1'}+\ldots}
It is easy to check that there are no other spurious
singularities in $x$ \foot{In principle one would expect the appearance of
singularities also in the branches in which eq. \thrbraone\ is satisfied
due to the symmetry between the variables $x_2$ and $x_3$ in \gfsaeq.
However, we remember that this symmetry has been explicitly broken
by he way in which the kernel $k(x,y)$ has been constructed.}.
Thus $K(x,y)$ is a differential of the third kind
on $\Sigma_4$. Taking into account also the simple pole in
$x=y$, the sum of all its residues vanishes as it should be on a
compact surface.

To study the singularities
with respect to the variables $y$, it is convenient to rewrite $K(x,y)$
in a slightly different form,
obtained by expanding in \gfker\ $f(y_1,y_2,x_3)$ and
$f(x_1,y_2,x_3)$ in powers of $x_3$ and $y_2$ respectively.
Since $f(x)$ is a polynomial in its arguments of degree three,
the expansions below:
\eqn\ggexpone{f(y_1,y_2,x_3)=\sum_{n=1}^3
{\partial^nf(y)\over\partial
y_3^n}{(x_3-y_3)^n\over n!}}
\eqn\ggexptwo{f(x_1,y_2,x_3)=\sum_{n=1}^3
{\partial^nf(x)\over\partial
x_2^n}{(x_2-y_2)^n\over n!}}
are exact.
As a consequence, inserting \ggexpone\ and \ggexptwo\
in \gfker, we have that
\eqn\gfserexpfor{K(x,y)={P(x,y)\over J^1(x)(x_1-y_1)}dx_1}
where
\eqn\pixy{
P(x,y)=\sum_{n=1}^3\left[
y_2 {\partial^nf(x)\over\partial
x_2^n}{(x_2-y_2)^{n-1}\over n!}-
x_3{\partial^nf(y)\over\partial
y_3^n}{(x_3-y_3)^{n-1}\over n!}\right]}
Now we exploit the fact
that the spurious divergences of $K(x,y)$ are located at the points in which
$y_1=\infty$ as the previous analysis has shown. Therefore, it
is convenient to keep in the kernel
only the contributions which diverge in $y_1=\infty$.
Using the formula ${1\over x_1-y_1}=-{1\over y_1}\sum\limits_{n=0}^\infty
\left({x_1\over y_1}\right)^n$, we find:
\eqn\expgf{K(x,y)\sim
{dx_1\over J^1(x)}\left[
{\cal A}_1(y)+x_1{\cal A}_2(y)+x_2{\cal A}_3(y)+x_3{\cal A}_4(y)\right]}
where
\eqn\calay{{\cal A}_1(y)=-y_1^{-1}\left(a_{03}^{(3)}y_2^3+
a_{02}^{(3)}y_2^2\right)}
\eqn\calby{{\cal A}_2(y)=-a_{03}^{(3)}{y_2^3\over y_1^2}-
a_{12}^{(3)}{y_2^2\over y_1}}
\eqn\calcy{{\cal A}_3(y)=-a_{03}^{(3)}{y_2^2\over y_1}}
\eqn\caldy{{\cal A}_4(y)=y_1^{-1}\left[
-a_{02}^{(2)}y_2^2+7y_3^2+3y_3h_1(y_1,y_2)+h_2(y_1,y_2)\right]}
We notice that $K(x,y)$ has the following behavior near the
spurious poles in $y_1=\infty$:
\eqn\asymbehav{K(x,y)\sim\sum_{i=1}^4\omega_i(x){\cal A}_i(y)}
where the $\omega_i(x)$ are holomorphic differentials.
As a matter of fact, using eqs. \thrbraone\ and \thrbratwo\ it is easy
to check that a basis of holomorphic differentials on $\Sigma_4$ is:
\eqn\holdifot{\omega_1(x)={dx_1\over J^1(x)}\qquad\qquad
\omega_2(x)={x_1dx_1\over J^1(x)}}
\eqn\holdiftf{\omega_3(x)={x_2dx_1\over J^1(x)}\qquad\qquad
\omega_4(x)={x_3dx_1\over J^1(x)}}
As a result, it is possible to conclude that the divergent part of
$K(x,y)$ given
in eq. \expgf\ is proportional to holomorphic differentials in $x$.
This property of the kernel \gfker\ will be crucial in subtracting the spurious
divergences from the amplitudes of the conformal field theories which will be
the subject of the next Section.
\newsec{FREE CONFORMAL FIELD THEORIES ON A GENERAL NON-HYPERELLIPTIC
CURVE OF GENUS FOUR}
In this Section we apply the previous results to the computation
of the amplitudes of the free conformal field theories appearing
in the action of bosonic strings on
a general non-hyperelliptic surface of genus four $\Sigma_4$.
These systems are well known and represent a good way to test
the generalized Weierstrass differential constructed in Section IV.
Apart from the above mentioned works in the case of hyperelliptic
and non-hyperelliptic curves, they
have been studied
by various authors and different methods on Riemann surfaces (see for instance
\ref\agr{
M. Bonini and R. Jengo, {\it Int. Jour. Mod. Phys.} {\bf A3}
(1988), 841;
E.
Date, M. Jimbo, M. Kashiwara and T. Miwa, In: Proc.
of International Symposium on Nonlinear Integrable Systems, Kyoto 1981,
M. Jimbo and T. Miwa (eds.), Singapore (1983); \prl{S. Saito}{36}{1987}{1819};
L. Alvarez-Gaum\'e, C. Gomez and C. Reina, New
Methods in String Theory, in: Superstrings '87, L. Alvarez-Gaum\'e
(ed.), Singapore, World Scientific 1988; \mpl{N. Ishibashi, Y. Matsuo and
Y. Ooguri}{2}{1987}{119}; N. Kawamoto, Y. Namikawa, A.
Tsuchiya and Y. Yamada, {\it Comm. Math. Phys.} {\bf 116} (1988),
247.}\nref\blmr{L. Bonora, A. Lugo,
M. Matone and J. Russo,
{\it Comm. Math. Phys.} {\bf 123} (1989), 329.}\nref\bontwo{\plb{L. Bonora, M.
Matone, F. Toppan and K. Wu}{224}{1989}{115}; \npbib{334}{1990}{717}; L. Bonora
and F. Toppan, {\it Rev. Math. Phys.} {\bf 4} (1992), 429.}
\nref\cvafa{\plb{C. Vafa}{190}{1987}{47}.}\nref\raina{A. K. Raina,
{\it Comm. Math. Phys.} {\bf 122} (1989), 625; {\it ibid.} {\bf 140}
(1991), 373; {\it Lett. Math. Phys.} {\bf 19} (1990), 1;
{\it Expositiones Mathematicae} {\bf 8} (1990), 227; {\it Helvetica
Physica Acta} {\bf 63} (1990), 694.}\nref\pdvpz{
P. di Vecchia,
{\it Phys. Lett.} {\bf B248} (1990), 329; P. Di Vecchia, F. Pezzella,
M. Frau, K. Hornfleck, A. Lerda and A. Sciuto, {\it Nucl. Phys.}
{\bf B332} (1989), 317; {\it ibid.} {\bf B333} (1990), 635; A. Clarizia
and F. Pezzella, {\it Nucl. Phys.} {\bf B298} (1988), 636;
\plb{G. Cristofano, R. Musto, F. Nicodemi and R. Pettorino}{217}{1989}
{59}.}\nref\russo{\npb{A. Lugo and J. Russo}{322}{1989}{210};
\plb{J.Russo}{220}{1989}{104}.}--\ref\semi{A. M.
Semikhatov, {\it Phys. Lett.} {\bf B212} (1988), 357;
O. Lechtenfeld, {\it Phys.
Lett} {\bf B232} (1989) 193; \npb{U. Carow-Watamura and S. Watamura}{288}
{1987}{500}; \npbib{301}{1988}{132}; \npbib{302}{1988}{149}; \npbib{308}{1988}
{143}; U. Carow-Watamura, Z. F. Ezawa, K.
Harada, A. Tezuka and S. Watamura, {\it Phys. Lett.} {\it B227}
(1989), 73.}).

First of all, we discuss the case of fermionic $b-c$ systems with integer spin
$\lambda=1,2$
In isothermal coordinates $p,\bar p$ where the metric becomes
conformally flat, their action is given by:
\eqn\bcsystems{S_{bc}=\int_{\Sigma_4}d^2p b\bar\partial c}
where $d^2p=idp\wedge d\bar p$ and $\bar\partial=
{\partial\over\partial\bar p}$.
The fields $b$ are meromorphic
tensors on $\Sigma_4$ with $\lambda$ lower indices, while
the fields $c$ are characterized by $\lambda-1$ upper indices.
From eq. \bcsystems\ one obtains the following
equations of motion:
\eqn\bcclaeqs{\bar\partial b=\bar\partial c= 0}
We start with the case $\lambda=2$.
We denote with $\phi_\mu$, $\mu=1,\ldots,9$, the holomorphic quadratic
differentials
which represent a basis of non-trivial solutions of \bcclaeqs.
If $p_1\ldots p_m$ and $q_1\ldots q_n$ are points in $\Sigma_4$,
the nonvanishing correlation functions of the $b-c$ systems
\eqn\bcltcf{{\cal G}_2(p;q)=
\langle\prod_{\alpha=1}^mb(p_\alpha) \prod_{\beta=1}^n
c(q_\alpha)\rangle}
should satisfy the relation $m-n=3g-3=9$.

The zeros and poles of ${\cal G}_2(p;q)$ are
determined by the physical properties of the fields. To specify their
locations and
orders it is convenient to introduce the concept of {\it divisor}. 
Let $\Delta [T]$ denote the divisor of a given meromorphic
tensor $T(p)$ on
the Riemann surface. If $T$ has zeros at $p=p_i$ of order $\mu_i$,
$i=1,\ldots,m_{zeros}$ and poles of order $\nu_j$
at $p=q_j$, $j=1,\ldots,n_{poles}$, then $\Delta [T]$ can be written as
follows:
\eqn\gendiv{\Delta [T]=\sum_{i=1}^{m_{zeros}}\mu_ip_i-
\sum_{j=1}^{n_{poles}}\nu_jq_j}
The correlation function \bcltcf\ is a tensor with two lower indices
in each variable $p_\alpha$, $\alpha=1,\ldots,m$ and the following divisor:
\eqn\divgip{\Delta_{p_\alpha}
[{\cal G}_2]=\sum_{\alpha'=1\atop\alpha'\ne \alpha}^m
p_{\alpha'}-\sum_{\beta'=1\atop\beta'\ne \beta}^n q_{\beta'}}
With respect to $q_\beta$, instead,
${\cal G}_2(p;q)$ is a vector with an upper index
and divisor:
\eqn\divgip{\Delta_{q_\beta}[{\cal G}_2]=
\sum_{\beta'=1\atop\beta'\ne \beta}^n q_{\beta'}-
\sum_{\alpha'=1\atop\alpha'\ne \alpha}^m p_{\alpha'}
}
To construct ${\cal G}_2(p;q)$ explicitly, we
represent $\Sigma_4$ as a ramified covering $C_4$ associated to the system
of algebraic equations \gfsaeq\ treated in the previous Section.
We begin by noting that a
meromorphic tensor $T(p)$ with $\lambda$ indices on the Riemann
surface $\Sigma_4$ corresponds in $C_4$ to a tensor
$T(x(p))=T_{x_1\ldots x_1}(x(p))dx_1^\lambda$, where
$x(p)=x_1(p),x_2(p),x_3(p)$ and $x_1(p)\in${\bf P}$_1$. A pole or a zero of
$T(p)$ in $p=q$ corresponds to a pole or a zero of $T(x(p))$ of the
same order in $x(p)=y(q)$.

A basis of independent holomorphic quadratic differentials is given by:
\eqn\diffsot{\phi_1(x)=\left({dx_1\over J^1(x)}\right)^2
\qquad\qquad \phi_2(x)=x_1\left({dx_1\over J^1(x)}\right)^2
\qquad\qquad \phi_3(x)=x_2\left({dx_1\over J^1(x)}\right)^2}
\eqn\diffsfs{\phi_4(x)=x_3\left({dx_1\over J^1(x)}\right)^2
\qquad\qquad \phi_5(x)=\left({x_1dx_1\over J^1(x)}\right)^2
\qquad\qquad \phi_6(x)=\left({x_2dx_1\over J^1(x)}\right)^2}
\eqn\diffsfs{\phi_7(x)=\left({x_3dx_1\over J^1(x)}\right)^2
\qquad\qquad \phi_8(x)=x_1x_2\left({dx_1\over J^1(x)}\right)^2
\qquad\qquad \phi_9(x)=x_1x_3\left({dx_1\over J^1(x)}\right)^2}
Clearly, a quadratic differential of the kind
$\phi_{10}(x)=x_2x_3\left({dx_1\over J^1(x)}\right)^2$ would not be
independent due to eq. \ggf.

It is also easy to check that the following quadratic differential:
\eqn\ktquadif{K_2(x,y)=-\left({dx_1\over J^1(x)}\right)^2\left[
{x_3f(y_1,y_2,x_3)\over(x_1-y_1)(x_3-y_3)}+
{y_2f(x_1,y_2,x_3)\over(x_1-y_1)(x_2-y_2)}\right]}
has only a simple pole at the point $x_i=y_i$, $i=1,2,3$.
In fact, $K_2(x,y)$ is obtained multiplying together
the Weierstrass kernel \gfker\ and the holomorphic differential
$\omega_1(x)$ of eq. \holdiftf.
The zeros of the latter cancel exactly the poles of $K(x,y)$. Indeed, since
$J^1(x)={\partial f\over\partial x_2}-x_1{\partial f\over\partial x_3}$
and $f$ is a polynomial of degree three in $x_3$, we have that
$J^1(x)\sim-{1\over x^{\prime 3}}$ at infinity in the three
branches in which eq. \thrbratwo\ is satisfied.

At this point we are
ready to write the correlation functions of the $b-c$ systems with
$\lambda=2$ on $C_4$:
$${\cal G}_2(p;q)=$$
\eqn\gtexplic{{\rm det}\left|
\matrix{K_2(x(p_1),y(q_1))&\ldots& K_2(x(p_1),y(q_1))&\phi_1(x(p_1))&\ldots
&\phi_9(x(p_1))
\cr
\vdots&\ddots&\vdots&\vdots&\ddots&\vdots\cr
K_2(x(p_m),y(q_1))&\ldots& K_2(x(p_m),y(q_n))&\phi_1(x(p_m))&\ldots
&\phi_9(x(p_m))
\cr
}
\right|
}
Due to the properties of determinants, the right hand side of the above
equation has the desired simple zeros whenever
$x_i(p_\alpha)=x_i(p_{\alpha'})$,
$\alpha,\alpha'=1,\ldots,m$ and $y_i(q_\beta)=y_i(q_{\beta'})$,
$\beta,\beta'=1,\ldots,n$, $i=1,2,3$.
Moreover, all the poles of
${\cal G}_2(p;q)$ are simple and occur at the points
in which $x_i(p_\alpha)=y_i(q_\beta)$.

In principle, there could be also spurious poles due to the fact that
$K_2(x,y)$ diverges when the variable $y$ becomes very large.
However, from what it has been discussed in the previous Section,
it is easy to realize that $K_2(x,y)$ has the following behavior in
$y_1=\infty$:
\eqn\nomse{K_2(x,y)\sim\sum_{i=1}^4\omega_i(x){\cal A}_i(y)}
where the functions ${\cal A}_i$ have been defined in eq. \asymbehav\
and diverge when $y_1\longrightarrow\infty$.
Terms of the form \nomse\ consist in a linear combination of quadratic
differentials, which does not contribute in the determinant of eq. \gtexplic.

Let us now treat the case $\lambda=1$. On the Riemann surface $\Sigma_4$ the
nonvanishing correlations functions are given by:
\eqn\gocorrfunc{{\cal G}_1(p,q)=
\langle\prod_{\alpha=1}^mb(p_\alpha) \prod_{\beta=1}^n
c(q_\alpha)\rangle}
where $m-n=g-1=3$.
Again, the poles and zeros of the above correlator are determined by the
physical properties of the $b-c$ fields. The divisors of ${\cal G}_1(p,q)$
are  similar to those of ${\cal G}_2(p,q)$ with simple zeros whenever
$p_\alpha=p_{\alpha'}$ or $q_\beta=p_{\beta'}$ and simple poles
when $p_\alpha=q_\beta$, with $\alpha,\alpha'=1,\ldots,m$,
$\beta,\beta'=1,\ldots,n$.

The explicit construction of ${\cal G}_1(p;q)$ on the ramified covering $C_4$
goes as follows. First of all, we define the differential:
\eqn\tkdiff{\nu_{y(q)y(q')}(x(p))=K(x(p),y(q))-K(x(p),y(q'))}
This is a differential of the third kind in $x(p)$, with two simple poles in
$x(p)=y(q)$ and $x(p)=y(q')$ and residues $+1$ and $-1$ respectively.
The spurious poles of $K(x(p),y(q))$ in $x_1=\infty$ are canceled against
the analogous poles of $K(x(p)-y(q'))$.
At this point it is possible to write the expression of ${\cal G}_1(p;q)$:
$${\cal G}_1(p;q)=$$
\eqn\goddd{{\rm det}\left|\matrix{
\nu_{y(q_1)y(q_n)}(x(p_1))&\ldots&\nu_{y(q_{n-1})y(q_n)}(x(p_1))&
\omega_1(x(p_1))&\ldots&\omega_4(x(p_1))\cr
\vdots&\ddots&\vdots&\vdots&\ddots&\vdots\cr
\nu_{y(q_1)y(q_n)}(x(p_m))&\ldots&\nu_{y(q_{n-1})y(q_n)}(x(p_m))&
\omega_1(x(p_m))&\ldots&\omega_4(x(p_m))\cr
}\right|
}
As in the case $\lambda=2$, eq. \asymbehav\ implies the absence of spurious
poles in the $y$ variables in \goddd.

To conclude our list of conformal field theories which appear in bosonic
string theory, we treat the scalar fields with action:
\eqn\scalact{S=\int_{\Sigma_4}d^2p\partial X\bar\partial X}
All the correlation functions of the scalar fields are obtained once the
following correlator is known:
\eqn\scagrefun{G(p;q,q')=\langle\partial_p X(p,\bar p)\left[
X(q,\bar q)-X(q',\bar q')\right]\rangle dp}
It turns out that $G(p;q,q')$ is a {\it canonical differential
of the third kind} 
uniquely determined by the following properties:
\item{$a$)} $G(p;q,q')$ has only
two simple poles in $p=q$ and $p=q'$ with residues $+1$ and $-1$ respectively.
\medskip
\item{$b$)} The integral function
$\int G(p;q,q')$ has purely imaginary periods when transported
around the $2g=8$ non-trivial homology cycles
of $\Sigma_4$.
\smallskip\noindent
On the algebraic curve $C_4$ the Green function \scagrefun\ 
can be written as a vector field
$G(x(p);y(q),y(q'))$, where $x(p)=x_1(p),x_2(p),x_3(p)$ etc.
$G(x(p);y(q),y(q'))$
coincides to the third kind differential $\nu_{y(q)y(q')}(x(p))$ defined
in eq. \tkdiff\ up to a linear combination of holomorphic differentials,
which
is fixed by requirement $b$).
In practice, since
it is hard to deal with integrals over homology cycles in the
case of algebraic curves, it is convenient
to formulate this requirement in terms of surface integrals.
%To this purpose, let us put
%\eqn\posaaa{G(x(p);y(q),y(q'))=\left(\tau_{y(q)y(q')}\right)_{x_1(p)}
%dx_1(p)}
Indeed, it is possible to show that $b$) is satisfied if and only
if $G(x(p);y(q),y(q'))$ fulfills the following Riemann bilinear identities:
\eqn\zscual{\int_{C_4} G(x(p);y(q),y(q'))\wedge
\overline{\omega_i(x(p))}=0\qquad\qquad i=1,\ldots,4}
The surface integrals over $C_4$ in \zscual\ can be interpreted
as integrals in a three dimensional complex space concentrated in the
solutions of eqs. \gfsaeq\ (see Appendix A).
At this point we are able to write
the Green function $G(x(p);y(q),y(q'))$ in terms of the third kind
differential \tkdiff\ and of the holomorphic differentials \holdifot-\holdiftf:
\eqn\gpqq{G(x(p);y(q),y(q'))={\rm det}\left|
\matrix{\nu_{y(q)y(q')}(x(p))&\omega_1(x(p))&\ldots&\omega_4(x(p))\cr
&&&\cr
\int_{C_4}\nu_{y(q)y(q')}\wedge\overline{\omega_1}&\int_{C_4}
\omega_1\wedge\overline{\omega_1}&\ldots&\int_{C_4}
\omega_1\wedge\overline{\omega_4}\cr
\vdots&\vdots&\ddots&\ldots\cr
\int_{C_4}\nu_{y(q)y(q')}\wedge\overline{\omega_4}&\int_{C_4}
\omega_4\wedge\overline{\omega_1}&\ldots&\int_{C_4}
\omega_4\wedge\overline{\omega_4}\cr
}\right|}
It is easy to see that the above differential of the third kind satisfies
requirement $a$) and the relations \zscual, which are equivalent to $b$).

\newsec{Conclusions}

In Section IV an analogue of the Weierstrass kernel has been constructed
on non-plane algebraic curves associated to the vanishing of two polynomials
$f$ and $g$.
The freedom of adding linear combinations of
differentials which do not change the behavior
in $x=y$ has been exploited in order to get two different,
but equivalent,
versions of generalized Weierstrass kernels.
The first version $K_{sym}(x,y)$, given in \ksymm, is
symmetric with respect to the
variables $x,y$ and with respect
to the exchange of $f$ with $g$.
This is in agreement with the fact that neither the
coordinates nor the functions $f$ and $g$ play a special role in the
equations which define the curve.
The alternative kernel $K(x,y)$ of eq. \fker\
has the advantage to have a more compact
expression in comparison to $K_{sym}(x,y)$, but part of the 
symmetry under coordinate permutations is lost. 
If the branched cover of an algebraic curve is considered, the generalized
Weierstrass kernel $K(x,y)$ has the simple form \simwgwk.

Furthermore, it has been verified that both kernels $K_{sym}(x,y)$ and 
$K(x,y)$ are third kind differentials with
a simple pole in the desired point $x=y$ of
the curve, in agreement with property \propess,
which characterizes the analogues
of the Cauchy kernel on Riemann surfaces \zvero.
In Section IV it has also been proved that spurious singularities
may only occur at the points in which one or more of the components
of the coordinates $x$ or $y$ approach infinity.
In the absence of a general algorithm to treat these spurious singularities
like that developed in the case of
plane curves in \naszejdwa--\naszed,
the terms to be subtracted in order to get the desired correlation
functions should be derived separately for any
given polynomials $f$ and $g$.
The example of a general non-hyperelliptic curve of genus four
has been explicitly worked out and the amplitudes of bosonic string theory
have been computed (see eqs. \gtexplic, \goddd\ and \gpqq).

Finally, nothing has been said about non-hyperelliptic
Riemann surfaces  obtained from the intersection of $n-1$ hypersurfaces in
${\bf P}^n$ with $n>3$. However, it is clear from eqs. \ksymm\ and \fker\
what should be the structure of the generalized
Weierstrass kernel on these curves. The Jacobians $J^i(x)$ should be
replaced by analogous Jacobians containing derivatives of
the $n-1$ polynomials
$f_1,\ldots,f_{n-1}$
with respect to any possible $n-1-$ dimensional subsets of the coordinates.
The numerators $N^i(x,y)$ will contain
a sum of products of polynomials
$f_1f_2\ldots f_{n-1}$, in which the dependence on the variables
$x_1,\ldots,x_n$ and $y_1,\ldots,y_n$
is chosen in such a way that the spurious poles
in the denominator given by the factor $\prod_{i=1}^n(x_i-y_i)$
are canceled, so that only the desired singularity in $x_i=y_i$,
$i=1,\ldots,n$ remains.
\vskip 1cm
\centerline{\bf Acknowledgements}
\vskip 1cm
This work has been completed during a visit
at the Center for Theoretical
Physics at MIT funded by
a Senior CNR-NATO grant, which is gratefully
acknowledged.
I am indebted to R. Jackiw for the warm hospitality
at CTP. I wish also to thank J. Sobczyk and W. Urbanik for
fruitful discussions.

\appendix{A}{Surface integrals over $C_4$}
In this Appendix we show that the surface
integrals appearing in eqs. \zscual\ and \gpqq\ can be expressed as integrals
in a three dimensional complex space in the presence of
Dirac $\delta-$functions which impose the constraints \gfsaeq.
We suppose here that $x\in${\bf C}$^3$, but the result is valid also for
a compact curve, in which case {\bf C}$^3$ has to be replaced by
{\bf P}$_1^3$ and the non-flat metric of
{\bf P}$_1$ should be taken into account.

Let us consider a surface integral of the kind:
\eqn\gsiotk{I=\int_{C_4}\rho}
where $\rho(x(p))$ is a $(1,1)-$form. In components
$\rho(x(p))= \rho_{x_1\bar x_1}(x(p))dx_1\wedge d\bar x_1$
and
\eqn\gsicomp{I=\int_{C_4}d^2x_1(p)\rho_{x_1\bar x_1}
(x_1(p),x_2(p),x_3(p))}
Since $C_4$ is a curve associated to the system of algebraic equations
\gfsaeq,
it is possible to rewrite the integral \gsiotk\
as follows:
\eqn\bdmsa{\int_{C_4}d^2x_1(p)\rho_{x_1\bar x_1}
(x(p))=\int_{{\bf C}^3}d^6x\left|J\left(
\matrix{f&g\cr x_2&x_3\cr}\right)\right|^2\delta^{(4)}(f,g)\rho_{x_1\bar x_1}
(x(p))
}
where $d^6x$ is the volume element in {\bf C}$^3$,
\eqn\jacdet{J\left(
\matrix{f&g\cr x_2&x_3\cr}\right)={\rm det}
\left|
\matrix{{\partial f\over \partial x_2}&{\partial g\over\partial x_2}\cr
&\cr
{\partial f\over \partial x_3}&{\partial g\over\partial x_3}\cr}\right|}
and the Dirac $\delta-$function $\delta^{(4)}(f,g)$ has been
defined using the formulas of
\ref\gfsh{I. M. Gel'fand and G. E. Shilov,
{\it Generalized functions}, Vol. I, Academic Press, New York and London,
1964.}.
After performing the change of variables $x_2,x_3\rightarrow f,g$, this
distribution becomes an usual four-dimensional $\delta-$function:
\eqn\fddf{\delta^{(4)}(f,g)=\left({\partial_f\partial_{\bar  f}+
\partial_g\partial_{\bar  g}\over2\pi^2}\right){1\over |f|^2+|g|^2}}
Now let us apply to eq. \bdmsa\ the inverse transformation which
brings back to the old coordinates. This can be done using the relations:
\eqn\chatravar{\partial_f={1\over g_{x_2}f_{x_3}-g_{x_3}f_{x_2}}
\left(g_{x_2}\partial_{x_3}-g_{x_3}\partial_{x_2}\right)
\qquad\qquad\partial_f={1\over g_{x_2}f_{x_3}-g_{x_3}f_{x_2}}
\left(f_{x_3}\partial_{x_2}-f_{x_2}\partial_{x_3}\right)}
Substituting the result in eq. \bdmsa, we obtain an explicit expression
of $I$ in terms of three dimensional complex integrals:
$$\int_{C_4}d^2x_1(p)\rho_{x_1\bar x_1}
(x(p))$$
\eqn\frappa{
={1\over 2\pi^2}\int d^6x
\rho_{x_1\bar x_1}(x(p))\left[\left|\left(
g_{x_2}\partial_{x_3}-g_{x_3}\partial_{x_2}\right)\right|^2
+
\left|\left(f_{x_3}\partial_{x_2}-f_{x_2}\partial_{x_3}\right)\right|^2\right]
{1\over |f|^2+|g|^2}
}

\footatend\vfill\supereject\immediate\closeout\rfile\writestoppt
%originariamente \baselineskip=14pt
\baselineskip=24pt\centerline{{\bf References}}\bigskip{\frenchspacing%
\parindent=20pt\escapechar=` \input refs.tmp\vfill\eject}\nonfrenchspacing
\bye